\documentclass[11pt]{amsart}
\usepackage{amsmath,amssymb,amsthm}
\usepackage[latin1]{inputenc}

\theoremstyle{plain}
\newtheorem{theorem}{Theorem}
\newtheorem{proposition}[theorem]{Proposition}
\newtheorem{lemma}[theorem]{Lemma}
\newtheorem{corollary}[theorem]{Corollary}

\theoremstyle{definition}

\newtheorem{remark}[theorem]{Remark}

\newtheorem*{theorem*}{Theorem}

\newcommand{\cb}{{\mathcal B}}

\newcommand{\ce}{{\mathcal E}}
\newcommand{\cf}{{\mathcal F}}
\newcommand{\cg}{{\mathcal G}}
\newcommand{\ch}{{\mathcal H}}
\newcommand{\ci}{{\mathcal I}}

\newcommand{\cl}{{\mathcal L}}

\newcommand{\cs}{{\mathcal S}}
\newcommand{\ct}{{\mathcal T}}
\newcommand{\cu}{{\mathcal U}}

\newcommand{\E}{{\mathbb E}}
\newcommand{\N}{{\mathbb N}}
\renewcommand{\P}{{\mathbb P}}

\newcommand{\R}{{\mathbb R}}
\newcommand{\Z}{{\mathbb Z}}

\newcommand{\ind}{{\bf 1}}

\newcommand{\norm}[1]{\mathop{\left\| #1 \right\|}\nolimits}
\newcommand{\val}[1]{\mathop{\left| #1 \right|}\nolimits}
\newcommand{\inv}[1]{\mathop{\frac{1}{ #1}}\nolimits}

\newcommand{\eps}{\varepsilon}
\newcommand{\vphi}{\varphi}

\newcommand{\cpol}{C_{\textup {\tiny pol}}^{\infty}(\R^{d})}
\newcommand{\cbo}{C_b^{\infty}(\R^{d})}
\newcommand{\clpol}{\cl\left(\cpol\right)}
\newcommand{\clbpol}{\cl_b\left(\cpol\right)}
\newcommand{\ntn}{\lfloor nt\rfloor/n}
\newcommand{\nt}{\lfloor nt\rfloor}
\newcommand{\nsn}{\lfloor ns\rfloor/n}

\begin{document}
\title[Euler Scheme and Tempered Distributions]{Euler Scheme and
  Tempered Distributions}
\date{\today}
\author{Julien Guyon}
\address{ENPC-CERMICS, 6 av. Blaise Pascal, Champs-sur-Marne, 77455
  Marne La Vall\'ee, France.} 
\email{julien.guyon@cermics.enpc.fr}

\begin{abstract}
Given a smooth $\mathbb{R}^d$-valued diffusion $(X^x_t,t\in[0,1])$ starting at point $x$, we study how fast the Euler scheme $X^{n,x}_1$
with time step $1/n$ converges in law to the random variable
$X^x_1$. Precisely, we look for which class of test functions $f$ the approximate expectation
$\E\left[f\left(X_1^{n,x}\right)\right]$ converges with speed $1/n$ to
$\E\left[f\left(X^x_1\right)\right]$.

When $f$ is smooth with polynomially growing derivatives
or, under a uniform hypoellipticity condition for
$X$, when $f$ is only measurable and bounded,
it is known that there exists a constant $C_1f(x)$ such that
\begin{equation}\label{eq:dl}
\E\left[f\left(X^{n,x}_1\right)\right] - \E\left[f\left(X^x_1\right)\right] =
C_1f(x)/n + O\left(1/n^2 \right).
\end{equation}

If $X$ is uniformly elliptic, we expand this result to the case when $f$
is a tempered distribution. In such a case, $\E\left[f\left(X^x_1\right)\right]$ (resp. $\E\left[f\left(X_1^{n,x}\right)\right]$) has to be understood as
  $\langle f,p(1,x,\cdot)\rangle$ (resp. $\langle f,p_n(1,x,\cdot)\rangle$) where $p(t,x,\cdot)$
  (resp. $p_n(t,x,\cdot)$) is the density of $X^x_t$ (resp. $X^{n,x}_t$). In
  particular, (\ref{eq:dl}) is valid when $f$ is a measurable function
  with polynomial growth, a Dirac mass or any derivative of a
  Dirac mass. We even show that (\ref{eq:dl}) remains valid when $f$ is a measurable
  function with exponential growth. Actually our results are symmetric
  in the two space variables $x$ and $y$ of the transition density and we prove that
\begin{equation*}
\partial^{\alpha}_x\partial^{\beta}_yp_n(t,x,y) -
\partial^{\alpha}_x\partial^{\beta}_yp(t,x,y) = \partial^{\alpha}_x\partial^{\beta}_y\pi(t,x,y)/n + r_n(t,x,y)
\end{equation*}
for a function $\partial^{\alpha}_x\partial^{\beta}_y\pi$ and a
$O(1/n^2)$ remainder $r_n$ which are shown to
have gaussian tails and whose dependence on $t$ is precised. We give applications to option pricing and hedging, proving
  numerical convergence rates for prices, deltas and gammas.
\end{abstract}


\maketitle

\section{Introduction and results}

Let $d,r\ge 1$ be two integers. Let $(\Omega,\cf,\P)$ be a probability
space on which lives a $r$-dimensional Brownian
motion $B$. We denote by $\cf_{t}=\sigma(B_{s},0\le
s\le t)$ the filtration generated by $B$. Let us give two functions
$b : \R^{d}\rightarrow\R^{d}$ and $\sigma :
\R^{d}\rightarrow\R^{d\times r}$. We systematically use (column) vector
and matrix notations, so that $b(x)$ should be thought of as a vector of
size $d$ and $\sigma(x)$ as a matrix of size $d\times r$. We denote
transposition by a star and define a $d\times d$ matrix-valued function
by putting $a=\sigma\sigma^*$. For a multiindex $\alpha\in\N^d$,
$\val\alpha=\alpha_1+\cdots+\alpha_d$ is its length and
$\partial^{\alpha}$ is the differential operator
$\partial^{\val\alpha}/\partial x^{\alpha_1}_1\cdots \partial x^{\alpha_d}_d$.
Equipping $\R^d$ with the euclidian norm $\norm{\cdot}$, we denote by
\begin{itemize}
\item $\cpol$ the set of infinitely differentiable functions
  $f:\R^{d}\rightarrow\R$ with polynomially growing derivatives of any
  order, i.e. such that for all
  $\alpha\in\N^d$, there exists $c\ge 0$ and $q\in\N$ such that for all $x\in\R^{d}$,
\begin{equation}\label{eq:C_infty_pol_Rd}
|\partial^{\alpha}\!f(x)| \le c\left(1+\|x\|^{q}\right),
\end{equation}
\item $\cbo$  the set of infinitely differentiable functions
  $f:\R^{d}\rightarrow\R$ with bounded derivatives of any
  order, i.e. such that $\partial^{\alpha}\!f\in L^{\infty}(\R^d)$ for all
  $\alpha\in\N^d$.
\end{itemize}

We shall make use of the following assumptions:
\begin{itemize}
\item[\bf{(A)}] For all $i\in\{1,\ldots,d\}$ and $j\in\{1,\ldots,r\}$, $b_i$
  and $\sigma_{i,j}$ belong to $C_{\textup {\tiny
      pol}}^{\infty}(\R^{d})$ and have bounded first derivatives.
\item[\bf{(B)}] For all $i\in\{1,\ldots,d\}$ and $j\in\{1,\ldots,r\}$, $b_i$
  and $\sigma_{i,j}$ belong to $C_b^{\infty}(\R^{d})$.
\item[\bf{(C)}] There exists $\eta>0$ such that for all $x,\xi\in\R^d$, $\xi^*a(x)\xi \ge \eta \|\xi\|^2$.
\end{itemize}
(C) is known as the uniform ellipticity condition.

It is well known that, given $x\in\R$, the hypothesis (A) guarantees the existence and the $\P$-almost
sure uniqueness of a solution $X^x=(X^x_t,t\ge 0)$ of the stochastic differential
equation (SDE)
\begin{equation}\label{eq:EDSmat}
X^x_t=x+\int_{0}^{t}b(X^x_{s})\:ds+\int_{0}^{t}\sigma(X^x_{s})\:dB_{s}.
\end{equation}

\subsection{Motivation}

Let us fix a time horizon $T>0$. Without loss of generality, we can and do assume that $T=1$. We try to estimate the law of
$X_{1}^x$. To do so, the most natural idea is to approach $X^x$ by its
Euler scheme of order $n\ge 1$, say $X^{n,x}=(X^{n,x}_t,t\ge 0)$, defined as follows. We
consider the regular subdivision
$\frak{S}_n=\{0=t^{n}_{0}<t^{n}_{1}<\cdots<t^{n}_{n-1}<t^{n}_{n}=1\}$ of
the interval $[0,1]$,
i.e. $t^{n}_{k}=k/n$, and we put $X^{n,x}_{0} = x$ and, for all $k\in\{0,\ldots,n-1\}$ and $t\in [t^{n}_{k},t^{n}_{k+1}]$,
\begin{eqnarray}\label{eq:def_euler}
X^{n,x}_{t} = X^{n,x}_{t^{n}_{k}}+b\left(X^{n,x}_{t^{n}_{k}}\right)(t-t^{n}_{k})
                            +\sigma\left(X^{n,x}_{t^{n}_{k}}\right)\left(B_{t}-B_{t^{n}_{k}}\right) .
\end{eqnarray}
Then the random variable $X_{1}^{n,x}$ is exactly simulatable and should
be close in law of $X_{1}^x$. Precisely, we measure the weak error between $X_{1}^{n,x}$ and $X_{1}^x$
by the quantities
\begin{equation*}
\Delta_1^nf(x)=\E\left[f\left(X^{n,x}_1\right)\right] - \E\left[f\left(X^x_1\right)\right]
\end{equation*}
and we try to find the largest space of test functions $f$ for which, for each $x$,
there exists a constant $C_1f(x)$ such that 
\begin{equation}\label{eq:dl2}
\Delta_1^nf(x) = C_1f(x)/n + O\left(1/n^2 \right).
\end{equation}

Practical interest of such an expansion has to be underlined (see, for
instance, \cite{jacod-kurtz-meleard-protter,talay-tubaro}). When (\ref{eq:dl2}) holds, one can use the Euler scheme plus a Monte-Carlo method to
estimate $\E\left[f\left(X^{x}_1\right)\right]$ and then, in a time
of order $nN$, gets an error of order $1/\sqrt{N} + 1/n$, where $N$
stands for the number of independants copies of $X^{n,x}_1$ generated by the Monte-Carlo
procedure. Given a tolerance $\eps\ll 1$, in order to minimize the time of
calculus, one should then choose $N=O\left(n^2\right)$ and gets a result
in a time of order $1/\eps^3$.

One can even do better using Romberg's extrapolation technique: if one
runs $N$ independant copies $(X^{2n,x}_{i,1},X^{n,x}_{i,1})$ of the couple $(X^{2n,x}_1,X^{n,x}_1)$,
which still requires a time of order $nN$, then computing $\inv
N\sum_{i=1}^N
(2f(X^{2n,x}_{i,1})-f(X^{n,x}_{i,1}))$
one gets an estimate of $\E\left[f\left(X^{x}_1\right)\right]$ whose
accuracy is of order $1/\sqrt{N} + 1/n^2$, since (\ref{eq:dl2}) implies
that
$\E[2f(X^{2n,x}_{1})-f(X^{n,x}_{1})]=\E[f(X^{x}_1)]+O(1/n^2)$.
Given a tolerance $\eps\ll 1$, one should now choose $N=O\left(n^4\right)$ and gets a result
in a time of order $1/\eps^{5/2}$.

\subsection{Previous results}

Using Itô expansions, D. \textsc{Talay} and L. \textsc{Tubaro}
\cite{talay-tubaro} have shown that (\ref{eq:dl2}) holds when $f\in
C_{\textup {\tiny pol}}^{\infty}(\R^{d})$ under condition
\begin{itemize}
\item[\bf{(B')}] The $b_i$'s and the $\sigma_{i,j}$'s are infinitely differentiable functions
  with bounded derivatives of any order $\ge 1$.
\end{itemize}
Hypothesis (B') is almost (B) but in (B') the functions $b_i$ and
$\sigma_{i,j}$ are not supposed bounded themselves. Using Malliavin calculus,
V. \textsc{Bally} and D. \textsc{Talay} \cite{bally-talay1} have
extended this result to the case of measurable and bounded $f$'s, with
the extra hypothesis that $X$ is uniformly hypoelliptic.
If (C) holds, $X_{1}^{n,x}$ and $X_{1}^x$ have densities, say $p_n(1,x,\cdot)$ and $p(1,x,\cdot)$ respectively (in this paper, densities are always taken with
respect to the Lebesgue measure). Then, for each pair $(x,y)$, the authors
\cite{bally-talay2} get an expansion of the error on the density itself
of the form
\begin{equation}\label{eq:dl_bally-talay2}
p_n(1,x,y)-p(1,x,y)=\pi(1,x,y)/n + \pi_n(1,x,y)/n^2.
\end{equation}
They also show that the principal error term $\pi$ and the remainder
$\pi_n$ have gaussian tails. Namely, they find constants $c_1\ge 0$ and $c_2>0$ such that for all $n\ge
1$ and $x,y\in\R^d$, $\val{\pi(1,x,y)} + \val{\pi_n(1,x,y)}\le
c_1\exp(-c_2\norm{x-y}^2)$. 

Besides, V. \textsc{Konakov} and E. \textsc{Mammen}
\cite{konakov-mammen} have proposed an analytical approach for this
problem based on the so-called parametrix method. If (B) and (C) hold, for each pair $(x,y)$, they get
an expansion of arbitrary order
$j$ of $p_n(1,x,y)$ but whose terms depend on $n$:
\begin{equation}\label{eq:dl_konakov-mammen}
p_n(1,x,y)-p(1,x,y)=\sum_{i=1}^{j-1}\pi_{n,i}(1,x,y)/n^i + O\left(1/n^{j}\right).
\end{equation}
They also prove that the coefficients have gaussian tails, uniformly in $n$: for each $i$,
they find constants $c_1\ge 0$ and $c_2>0$ such that for all $n\ge 1$
and $x,y\in\R^d$, $\val{\pi_{n,i}(1,x,y)}\le
c_1\exp(-c_2\norm{x-y}^2)$. To do so, the authors use upper bounds on
the partial derivatives of $p$ - which they find in \cite{friedman1} - and prove analogous bounds
on $p_n$'s ones.

A link with generalized Watanabe distributions on Wiener's space is
exhibited in \cite{malliavin-thalmaier}. For the general case of Lévy driven stochastic differential equations,
(\ref{eq:dl2}) holds under regularity assumptions on $f$ and integrability
conditions on the Lévy process, see
\cite{jacod-kurtz-meleard-protter,protter-talay}. The rate of
convergence of the process $(X^{n,x}_t-X^x_t,t\in[0,1])$ is given in
\cite{jacod-protter,jacod}. As for the simulation of densities, see for
instance \cite{kohatsu-pettersson}.

\subsection{Purpose and method}\label{sec:purpose}

Equations (\ref{eq:dl_bally-talay2}) and (\ref{eq:dl_konakov-mammen})
can be seen as expansions of
$\Delta_1^nf(x)=\E\left[f\left(X^{n,x}_1\right)\right] -
\E\left[f\left(X^x_1\right)\right]$ in the special case when
$f=\delta_y$, the Dirac mass at point $y\in\R^d$. We aim at giving a
precise sense to such quantities when $f$ is any tempered distribution,
and at proving that expansions in powers of $1/n$ remain valid in this extremely general
setting. Moreover, we will derive expansions that are valid not only for
$t=1$, but also for any time $t\in(0,1]$, the stepsize $1/n$ being
fixed, and we shall make explicit, in these expansions, the way the
coefficients and the remainders depend on $t$, $f$
and $x$.

To get these precise results, we shall place ourselves in a strong
situation. Namely, we will assume infinite regularity and boundedness of
the coefficients of the SDE (\ref{eq:EDSmat}), that is condition (B),
and uniform ellipticity, that is condition (C). The reason for this is the
following. Let us write $P_tf(x)=\E[f(X^x_t)]$ and
$P^n_tf(x)=\E\left[f(X^{n,x}_t)\right]$. We first expand
$\Delta^n_t=P^n_t-P_t$ as an endomorphism of $\cpol$, in powers of $1/n$. This can be done under
nothing more than hypothesis (A), see Theorems
\ref{thm:principal_section_esp} and \ref{thm:principal_section_esp_km}
in Section \ref{sec:preliminary}. The coefficients in
these expansions are operators of the form $\int_0^tP_sDP_{t-s}\: ds$ or
$\sum_{t^n_k<t}P^n_{t^n_k}DP_{t-t^n_k}$, where $D$ is a differential
operator. Now, under (B) and (C), both $X_{t}^{x}$ and $X_{t}^{n,x}$
have regular densities, say $p(t,x,\cdot)$ and $p_n(t,x,\cdot)$, with
gaussian tails, as soon as
$t>0$, so that we may express these operators as integral operators on $\R^d$. For
instance, for $f\in\cpol$, $x\in\R^d$ and $0<s<t$,
\begin{eqnarray*}
P_sDP_{t-s}f(x) &=& \int_{\R^d}p(s,x,z)DP_{t-s}f(z)\: dz \\
         &=& \int_{\R^d}
f(y)\left(\int_{\R^d}p(s,x,z)D(p(t-s,\cdot,y))(z)\:dz\right)\: dy.
\end{eqnarray*}
Now the expansions read on the density itself, with coefficients of the
form
\begin{eqnarray}
\pi(t,x,y) &=& \int_0^t\int_{\R^d}p(s,x,z)D(p(t-s,\cdot,y))(z)\: dzds, \label{eq:coeff1}\\
\mbox{or }\quad\pi_n(t,x,y) &=& \sum_{t^n_k<t}\int_{\R^d}p_n(t^n_k,x,z)D(p(t-t^n_k,\cdot,y))(z)\: dz.\label{eq:coeff2}
\end{eqnarray}

At this step, the key point is to prove that these coefficients, as well
as any of their spatial derivatives, have gaussian tails (see
Proposition \ref{prop:convol_esp}). To do so,
we split the above time integral (resp. sum) depending on whether $s$
(resp. $t^n_k$) is small or large, and integrate by parts in the latter
case. This is very similar to V. \textsc{Bally} and D. \textsc{Talay}'s
technique \cite{bally-talay1}, but they use the Malliavin calculus
integration by parts formula whereas we only use the genuine one. Then we use upper bounds on
the partial derivatives of $p$ and $p_n$, as is done in
V. \textsc{Konakov} and E. \textsc{Mammen}'s work
\cite{konakov-mammen}. Here the uniform ellipticity hypothesis is
crucial: it provides upper bounds that have enough quality in $t$ to allow us to
conclude.

The same analysis, with a bit more work, can be done for the remainders. We then
get \emph{functional} expansions of the form
\begin{equation}\label{eq:functional_exp}
p_n - p = \pi/n + \pi_n/n^2 \qquad  \mbox{ or }\qquad p_n-p=\sum_{i=1}^{j}\pi_{n,i}/n^i
\end{equation}
where $\pi$,
$\pi_n$ and the $\pi_{n,i}$'s and all their spatial derivatives have gaussian tails, uniformly in
$n$. We then achieve to give a distributional sense to expansion (\ref{eq:dl2}) by a duality approach: any tempered distribution can be
integrated or bracketed in the variable $y$ with the
expansions. Theorems \ref{thm:principal}, \ref{thm:principal2} and \ref{thm:principal_km} provide
precise statements, see Section \ref{sec:main}. 

\subsection{A first series of results}\label{sec:firstseries}

Stating Theorems \ref{thm:principal} and \ref{thm:principal_km} requires
a bit of preparation, namely defining appropriate functional spaces in
which will live the coefficients $\pi$, $\pi_n$ and $\pi_{n,i}$ in
expansions (\ref{eq:functional_exp}). Before doing
this, to encourage the reader, we would like to state a series of easy consequences of
Theorem \ref{thm:principal}, including an application to financial markets. They will be proved in Section
\ref{sec:main}. The function $\pi$ which appears in them is the principal
{\it functional} error term. It is defined by (\ref{eq:pi})-(\ref{eq:L2star}). Note that analogous corollaries can be derived from
Theorem \ref{thm:principal_km} as well. The first result gives the rate
of convergence of the spatial derivatives of the density:

\begin{proposition}\label{cor:principal1}
Under (B) and (C), for all $\alpha,\beta\in\N^d$, there exists $c_1\ge 0$ and $c_2>0$ such that for all $n\ge 1$, $t\in(0,1]$ and $x,y\in\R^{d}$,
\begin{equation*}
\partial^{\alpha}_x\partial^{\beta}_yp_n(t,x,y) - \partial^{\alpha}_x\partial^{\beta}_yp(t,x,y) =
\inv n \partial^{\alpha}_x\partial^{\beta}_y\pi(t,x,y) + r_n(t,x,y)
\end{equation*}
and
\begin{equation*}
\val{r_n(t,x,y)}\le c_1 n^{-2}t^{-(\val\alpha+\val\beta+d+4)/2}\exp\left(-c_2\norm{x-y}^2/t\right).
\end{equation*}
\end{proposition}

The next proposition states that (\ref{eq:dl2}) is valid for measurable
and polynomially growing $f$'s:

\begin{proposition}\label{cor:f_mes_pol}
Assume (B) and (C). Let $f:\R^d\rightarrow\R$ be a measurable function
such that there exists $c'\ge 0$ and $q\in\N$ such that for all
$x\in\R^d$, $\val{f(x)}\le c'(1+\norm{x}^q)$. Then there exists $c\ge 0$ such that for all $n\ge 1$, $t\in(0,1]$ and $x\in\R^d$,
\begin{equation}\label{eq:dl_croiss_pol}
\E[f(X^{n,x}_t)] - \E[f(X^x_t)] =
\inv n \int_{\R^d}f(y)\pi(t,x,y)\: dy+ r_n(t,x)
\end{equation}
and
\begin{equation*}
\val{r_n(t,x)}\le
cn^{-2}t^{-2}\left(1+\norm{x}^q\right).
\end{equation*}
\end{proposition}

As far as extending the class of $f$'s for which (\ref{eq:dl2}) holds is
concerned, we can even do better. Indeed, if for $\mu\in(0,2)$ we denote by $\ce_{\mu}$ the set
of all measurable functions $f:\R^d\rightarrow\R$ such that there exists
$c_1,c_2\ge0$ such that for all $y\in\R^d$,
\begin{equation*}
\val{f(y)} \le c_1 \exp\left(c_2\norm y ^{\mu}\right),
\end{equation*}
we have 
\begin{proposition}\label{cor:exp}
Under (B) and (C), for all $\mu\in(0,2)$ and $f\in\ce_{\mu}$, there exists $c_1,c_2\ge 0$ such that for all $n\ge 1$, $t\in(0,1]$ and $x\in\R^d$, $f(X^x_t)$ and $f(X^{n,x}_t)$ are integrable and
\begin{equation}\label{eq:dl_ce_mu}
\E[f(X^{n,x}_t)] - \E[f(X^x_t)] = \inv n \int_{\R^d}f(y)\pi(t,x,y)\: dy + r_n(t,x)
\end{equation}
with
\begin{equation*}
\val{r_n(t,x)}\le
c_1n^{-2}t^{-2}\exp\left(c_2\norm x ^{\mu}\right).
\end{equation*}
\end{proposition}

In particular, (\ref{eq:dl2}) remains true under (B) and (C) when $f\in\ce=\cup_{\mu\in(0,2)}\ce_{\mu}$.
More generally, Theorem \ref{thm:principal} leads to
\begin{proposition}\label{cor:d/dx_exp}
Under (B) and (C), for all $\alpha\in\N^d$, $\mu\in(0,2)$ and $f\in\ce_{\mu}$, there exists $c_1,c_2\ge 0$ such that for all $n\ge 1$, $t\in(0,1]$ and $x\in\R^d$,
\begin{equation}\label{eq:dl_ce_mu_derive}
\partial^{\alpha}_x\E[f(X^{n,x}_t)] - \partial^{\alpha}_x\E[f(X^x_t)] = \inv n \int_{\R^d}f(y)\partial^{\alpha}_x\pi(t,x,y)\: dy + r_n(t,x)
\end{equation}
with
\begin{equation*}
\val{r_n(t,x)}\le
c_1n^{-2}t^{-(\val\alpha+4)/2}\exp\left(c_2\norm x ^{\mu}\right).
\end{equation*}
\end{proposition}

This result can now be used in the context of financial markets.

\subsection{Application to option pricing and hedging}\label{sec:finance}

Let $S^v=(S^{v,1},\ldots,S^{v,d})$ be a basket of assets satisfying
\begin{equation*}
\frac{dS^{v,i}_t}{S^{v,i}_t} = \mu_i(S^v_t)\: dt +
\sum_{j=1}^{r}\sigma_{i,j}(S^v_{t})\: dB^{j}_{t}, \qquad S^{v,i}_0=v^i>0,
\end{equation*}
with $\mu,\sigma\in\cbo$ and $\sigma$ satisfying (C). Given a measurable and polynomially
growing function $\phi$, we try to estimate the price $\mbox{Price}=\E[\phi(S^v_t)]$, the deltas
$\mbox{Delta}_i=\partial^{e_i}_v\E[\phi(S^v_t)]$ and the gammas
$\mbox{Gamma}_{i,j}=\partial^{e_i+e_j}_v\E[\phi(S^v_t)]$ of the european option
of maturity $t$ and payoff $\phi$ ($(e_1,\ldots,e_d)$ is the canonical base of $\R^d$). To do so, let us set $x=\ln v$ (i.e. $x^i=\ln v^i$) and
$X^{x,i}_t=\ln(S^{v,i}_t)$. Then $X$ is the solution of (\ref{eq:EDSmat})
with $b=\mu-\norm{\sigma}^2/2\in\cbo$, where
$\norm{\sigma}^2_i(x)=\sum_{j=1}^{r}\sigma_{i,j}^2(x)$.
If we set $\exp(x)=(\exp(x^1),\ldots,\exp(x^d))$ and $f(x)=\phi(\exp(x))$, we
define a function $f\in\ce_1$ and, since $\mbox{Price}=\E[f(X^x_t)]$, (\ref{eq:dl_ce_mu}) leads to
\begin{equation*}
\mbox{Price}^n-\mbox{Price} = C^{\mbox{\tiny Price}}_t\phi(v)/n + O\left(n^{-2}t^{-2}\exp\left(c_2\norm{\ln v}\right)\right),
\end{equation*}
where $\mbox{Price}^n$ stands for the approximated price
$\E[f(X^{n,x}_t)]$ and
\begin{equation*}
C^{\mbox{\tiny Price}}_t\phi(v)=\int_{(\R_+^*)^d}\phi(u)\frac{\pi(t,\ln v,\ln u)}{u_1\cdots u_d}\: du.
\end{equation*}
Besides, if we set $\mbox{Delta}^n_i= \partial^{e_i}_v\E[f(X^{n,\ln v}_t)]$ and $\mbox{Gamma}^n_{i,j}= \partial^{e_i+e_j}_v\E[f(X^{n,\ln v}_t)]$, (\ref{eq:dl_ce_mu_derive}) shows that
\begin{eqnarray*}
\mbox{Delta}^n-\mbox{Delta} &=&  C^{\mbox{\tiny Delta}}_t\phi(v)/n + O\left(n^{-2}t^{-5/2}\exp\left(c_2\norm{\ln v}\right)\right), \\
\mbox{Gamma}^n-\mbox{Gamma} &=&  C^{\mbox{\tiny Gamma}}_t\phi(v)/n + O\left(n^{-2}t^{-3}\exp\left(c_2\norm{\ln v}\right)\right),
\end{eqnarray*}
where
\begin{eqnarray*}
C^{\mbox{\tiny Delta}}_t\phi(v)_i&=&\inv{v_i}\int_{(\R_+^*)^d}\phi(u)\frac{\partial^{e_i}_2\pi(t,\ln v,\ln u)}{u_1\cdots u_d}\: du, \\
C^{\mbox{\tiny Gamma}}_t\phi(v)_{i,j}&=&\inv{v_iv_j}\int_{(\R_+^*)^d}\phi(u)\frac{\partial^{e_i+e_j}_2\pi(t,\ln v,\ln u)-\ind_{\{i=j\}}\partial^{e_i}_2\pi(t,\ln v,\ln u)}{u_1\cdots u_d}\: du.
\end{eqnarray*}
Eventually we have proved that applying the Euler scheme of order $n$ to the
logarithm of the underlying leads to approximations of the price, the deltas
and the gammas which converge to the true price, deltas and gammas with
speed $1/n$, at least when the drift and volatility of the underlying satisfy (B)
and (C), which in the context of financial markets seems not to be a restricting hypothesis.
Note that the principal part of the error explodes as $t$ tends to 0 as
$t^{-1/2}$ for the prices, $t^{-1}$ for the deltas and $t^{-3/2}$ for
the gammas.

\subsection{Some functional spaces}\label{sec:def_spaces}

In order to state our main results (Proposition \ref{prop:convol_esp}
and Theorems \ref{thm:principal} and \ref{thm:principal_km}) precisely and shortly, let us introduce some families of functional
spaces. Functional expansions like (\ref{eq:functional_exp}) will take
place in such spaces. For $l\in\Z$, we first define $\cg_l(\R^{d})$ as the set of all measurable functions
$\pi:(0,1]\times\R^{d}\times\R^{d}\rightarrow\R$ such that
\begin{itemize}
\item for all $t\in (0,1]$, $\pi(t,\cdot,\cdot)$ is infinitely differentiable,
\item for all $\alpha,\beta\in \N^{d}$, there exists two constants $c_1\ge 0$ and $c_2>0$ such that for all $t\in(0,1]$ and $x,y\in\R^d$,
\begin{equation}\label{eq:def_cgl1}
\val{\partial^{\alpha}_x\partial^{\beta}_y\pi(t,x,y)} \le c_1 t^{-(\val\alpha+\val\beta+d+l)/2}\exp\left(-c_2\norm{x-y}^2/t\right).
\end{equation}
\end{itemize}
We say that a subset $\cb\subset\cg_l(\R^{d})$ is bounded if, in (\ref{eq:def_cgl1}), $c_1$ and $c_2$ can be chosen independently on $\pi\in\cb$. We also introduce the space $\cg(\R^{d})$ defined in the same way as $\cg_l(\R^{d})$ with (\ref{eq:def_cgl1}) replaced by the following two conditions:
\begin{eqnarray}
\val{\partial^{\alpha}_x\partial^{\beta}_y\pi(t,x,y)} &\le& c_1 t^{-(\val\alpha+\val\beta+d)/2}\exp\left(-c_2\norm{x-y}^2/t\right),\label{eq:def_cg1} \\
\val{\partial^{\alpha}_x\left(\pi\left(t,x,x+y\sqrt{t}\right)\right)} &\le& c_1t^{-d/2}\exp\left(-c_2\norm{y}^2\right).\label{eq:def_cg2}
\end{eqnarray}
Note that we may always take the couple of constants $(c_1,c_2)$ to be
the same in both equations (\ref{eq:def_cg1}) and
(\ref{eq:def_cg2}). Indeed, if they hold with two couples $(c'_1,c'_2)$
and $(c''_1,c''_2)$, they both hold with $(c_1,c_2)$ if we take
$c_1=c'_1\vee c''_1$ and $c_2=c'_2\wedge c''_2$. We say that a subset
$\cb\subset\cg(\R^{d})$ is bounded if, in (\ref{eq:def_cg1}) and
(\ref{eq:def_cg2}), $c_1$ and $c_2$ can be chosen independently on
$\pi\in\cb$. Note that in equation (\ref{eq:def_cg2}), the
upper bound keeps the same quality in $t$, namely $t^{-d/2}$, whatever the ``number'' $\alpha$ of times one differentiates
the mapping $x\mapsto \pi\left(t,x,x+y\sqrt{t}\right)$. This will be
crucial when proving Proposition \ref{prop:convol_esp}.

It is convenient to extend these definitions to mappings that also
depend on an intermediate time $s\in(0,t)$. To do so, let us denote by $\ct_1$ the
unit triangle $\{(s,t)\in\R^2|0<s<t\le 1\}$ and, for $l\in\Z$, let us define $\ch_l(\R^{d})$ as the space of measurable functions
$\rho:\ct_1\times\R^{d}\times\R^{d}\rightarrow\R$ such that
\begin{itemize}
\item for all $(s,t)\in\ct_1$, $\rho(s,t,\cdot,\cdot)$ is infinitely differentiable,
\item for all $\alpha,\beta\in \N^{d}$, there exists two constants $c_1\ge 0$ and $c_2>0$ such that for all $(s,t)\in\ct_1$ and $x,y\in\R^d$,
\begin{equation}\label{eq:def_chl1}
\val{\partial^{\alpha}_x\partial^{\beta}_y\rho(s,t,x,y)} \le c_1 t^{-\left(\val\alpha+\val\beta+d+l\right)/2}\exp\left(-c_2\norm{x-y}^2/t\right).
\end{equation}
\end{itemize}
Again we say that a subset $\cb\subset\ch_l(\R^{d})$ is bounded if, in
(\ref{eq:def_chl1}), $c_1$ and $c_2$ can be chosen independently on
$\rho\in\cb$. We also introduce the space $\ch(\R^{d})$ which is defined
in the same way as $\ch_l(\R^{d})$ with (\ref{eq:def_chl1}) replaced by
\begin{eqnarray}
\val{\partial^{\alpha}_x\partial^{\beta}_y\rho(s,t,x,y)} &\le& c_1 t^{-\left(\val\alpha+\val\beta+d\right)/2}\exp\left(-c_2\norm{x-y}^2/t\right),\label{eq:def_ch1} \\
\val{\partial^{\alpha}_x\left(\rho\left(s,t,x,x+y\sqrt{t}\right)\right)} &\le& c_1t^{-d/2}\exp\left(-c_2\norm{y}^2\right),\label{eq:def_ch2}
\end{eqnarray}
and we say that a subset $\cb\subset\ch(\R^{d})$ is bounded if, in
(\ref{eq:def_ch1}) and (\ref{eq:def_ch2}), $c_1$ and $c_2$ can be chosen independently on
$\rho\in\cb$. Again we may always choose the couple $(c_1,c_2)$ to be
the same in both equations (\ref{eq:def_ch1}) and
(\ref{eq:def_ch2}). Note that the upper bounds in (\ref{eq:def_chl1}), (\ref{eq:def_ch1}) and
(\ref{eq:def_ch2}) are exactly the same as the ones in (\ref{eq:def_cgl1}), (\ref{eq:def_cg1}) and
(\ref{eq:def_cg2}). In particular, they do not depend on $s$.

Eventually, for $\pi_1,\pi_2\in\cg(\R^{d})$, $g\in\cbo$ and $\gamma\in\N^d$, we define a function $\pi_1*_{g,\gamma}\pi_2$ on $\ct_1\times\R^{d}\times\R^{d}$ by putting
\begin{equation*}
  \left(\pi_1*_{g,\gamma}\pi_2\right)(s,t,x,y)=\int_{\R^{d}}g(z)\pi_1(s,x,z)\partial^{\gamma}_2\pi_2(t-s,z,y)\: dz.
\end{equation*}
Notation $\partial_2$ means differentiation with respect to the second
argument, here $z$. Operation $*_{g,\gamma}$ is a space convolution
which naturally appears when developping the differential operator $D$
in equations (\ref{eq:coeff1}) and (\ref{eq:coeff2}).

\subsection{Main results}\label{sec:main}

We are now able to state our main results as follows.

\begin{proposition}\label{prop:convol_esp}
Let $\cb_1$ and $\cb_2$ be two bounded subsets of $\cg(\R^{d})$,
$g\in\cbo$ and $\gamma\in\N^d$. Then
\begin{itemize}
\item[(i)] $\{\pi_1*_{g,\gamma}\pi_2|\pi_1\in\cb_1,\pi_2\in\cb_2\}$ is a
  bounded subset of $\ch_{\val{\gamma}}(\R^{d})$,
\item[(ii)] $\{\pi_1*_{g,0}\pi_2|\pi_1\in\cb_1,\pi_2\in\cb_2\}$ is a bounded subset of $\ch(\R^{d})$.
\end{itemize}
\end{proposition}


\begin{theorem}\label{thm:principal}
Under (B) and (C),
\begin{itemize}
\item[(i)] for all $t\in (0,1]$ and $x\in\R^d$, $X_{t}^{x}$ has a density $p(t,x,\cdot)$ and $p\in\cg(\R^{d})$,
\item[(ii)] for all $t\in (0,1]$, $x\in\R^d$ and $n\ge 1$, $X_{t}^{n,x}$ has a density $p_n(t,x,\cdot)$ and $(p_n,n\ge 1)$ is a bounded sequence in $\cg(\R^{d})$,
\item[(iii)] there exists $\pi\in\cg_1(\R^{d})$ and a bounded sequence $(\pi_n,n\ge 1)$ in $\cg_4(\R^{d})$ such that for all $n\ge 1$,
\begin{equation}\label{eq:principal}
p_n-p= \pi/n + \pi_n/n^{2}.
\end{equation}
\end{itemize}
\end{theorem}

These results are proved in Section \ref{sec:proof_thm_principal}. In
Theorem \ref{thm:principal}, statement (i) is already known, see
\cite{friedman1}, Theorem 7, page 260, and statement (ii) has
essentially been proved in \cite{konakov-mammen}. As explained in
Section \ref{sec:purpose}, Proposition \ref{prop:convol_esp}, together
with these two statements, is the key to derive statement (iii).

The function $\pi$ can be expressed in terms of $p$ by
\begin{equation}\label{eq:pi}
\pi(t,x,y) = \inv 2 \int_{0}^{t}\int_{\R^d}p(s,x,z)L^*_2(p(t-s,\cdot,y))(z)\: dzds,
\end{equation}
where the differential operator $L^*_2$ is explicitely given in terms of the functions $a$ and $b$ by
\begin{multline}\label{eq:L2star}
-L^*_2=\sum_{i=1}^d \left(b\cdot\nabla b_i + \inv 2 \mbox{tr}\left(a\nabla^2 b_i \right) \right)\partial_{i} \\+ \sum_{i,j=1}^d \left(\inv 2 b\cdot\nabla a_{i,j} + a_j\cdot\nabla b_i + \inv 4 \mbox{tr}\left(a\nabla^2 a_{i,j}\right) \right)\partial_{ij} + \inv 2 \sum_{i,j,k=1}^d a_k\cdot\nabla a_{i,j}\partial_{ijk}.
\end{multline}
Here, $\cdot$, $a_k$, tr, $\nabla$ and $\nabla^2$ respectively stand for
the inner product in $\R^d$, the $k$-th column of $a$, the trace of a
matrix, the gradient vector and the hessian matrix. In the case when
$t=1$, (\ref{eq:pi}) agrees with  V. \textsc{Bally} and
D. \textsc{Talay}'s expression for $\pi$ (\cite{bally-talay2},
definition 2.2, page 100), but seems preferable because it does not
involve differentiation with respect to $t$ and makes explicit that the
space differential operator $L_2^*$ is of order less than 3, when
V. \textsc{Bally} and D. \textsc{Talay}'s operator $\cu$ involves a
fourth order differentiation in space.

\vspace{.5cm}

We shall now prove that if $X$ is elliptic the expansion (\ref{eq:dl2}) is valid in the
very general case when $f$ is a tempered distribution. Let us denote by $\cs(\R^{d})$ Schwartz's space, i.e. the space of
infinitely differentiable functions $\vphi:\R^{d}\rightarrow\R$ such
that $x\mapsto x^{\alpha}\partial^{\beta}\vphi(x)\in L^{\infty}(\R^d)$
for all $\alpha,\beta\in\N^{d}$ ($x^{\alpha}$ stands for $x_1^{\alpha_1}\cdots
x_d^{\alpha_d}$),
and let us denote by $\cs'(\R^{d})$ the space of tempered
distributions. The seminorms $(N_q,q\in\N)$ are defined on
$\cs(\R^{d})$ by
 \begin{equation*}
N_q(\vphi)=\sum_{|\alpha|\le q,|\beta|\le
  q}\sup_{x\in\R^{d}}\val{x^{\alpha}\partial^{\beta}\vphi(x)},
\end{equation*}
and the order $\#S$ of $S\in\cs'(\R^d)$ is the smallest integer $q$ such that
there is a $c\ge 0$ such that $\val{\langle S,\vphi\rangle}\le
cN_q(\vphi)$ for all $\vphi\in\cs(\R^d)$.
Note that whenever $\pi\in\cg_l(\R^d)$, $\pi(t,x,\cdot)$ and $\pi(t,\cdot,y)$  belong to
$\cs(\R^{d})$. More precisely, for $\cb\subset\cg_l(\R^d)$ bounded, there
exists $c\ge0$ such that for all $\pi\in\cb$, $t\in(0,1]$ and $x,y\in\R^d$,
\begin{equation*}
N_q(\pi(t,x,\cdot))\le
ct^{-(d+l+q)/2}\left(1+\norm{x}^q\right) \!\qquad
\mbox{and} \qquad\! N_q(\pi(t,\cdot,y))\le
ct^{-(d+l+q)/2}\left(1+\norm{y}^q\right).
\end{equation*}

Applying a tempered distribution $S$
to (\ref{eq:principal}), $t$ and $x$ or $t$ and $y$ being fixed, we immediately deduce from Theorem \ref{thm:principal}

\begin{theorem}\label{thm:principal2}
Under (B) and (C), for all $S\in\cs'(\R^{d})$, there exists $c\ge 0$ such that for all $n\ge 1$, $t\in(0,1]$ and $x,y\in\R^d$,
\begin{eqnarray*}
\langle S,p_n(t,x,\cdot) \rangle - \langle S,p(t,x,\cdot) \rangle &=&
\inv n \langle S,\pi(t,x,\cdot) \rangle+ r_n'(t,x),\\
\langle S,p_n(t,\cdot,y) \rangle - \langle S,p(t,\cdot,y) \rangle &=&
\inv n \langle S,\pi(t,\cdot,y) \rangle + r_n''(t,y),
\end{eqnarray*}
and
\begin{equation*}
\val{r_n'(t,x)}+\val{r_n''(t,x)}\le
cn^{-2}t^{-(d+4+\#S)/2}\left(1+\norm{x}^{\#S}\right).
\end{equation*}
\end{theorem}

Let us define $\E\left[S(Y)\right]$ by $\langle S,p_Y
  \rangle$ when $S\in\cs'(\R^{d})$ and $Y$ is a random variable with
  density $p_Y\in\cs(\R^{d})$. Note that, when $S$ is a measurable and
  polynomially growing function, this definition coincides with the usual expectation. We then have proved that, under (B) and (C), (\ref{eq:dl2}) is
  valid for $f$'s being only tempered distributions, and not only for
  $t=1$, but also for any time $t\in(0,1]$, and we have even
  precised the way the $O(1/n^2)$ remainder depends on $t$, $f$ and $x$. Precisely, this remainder grows
  slower than $\norm x ^{\#f}$ as $x$ tends to infinity, and
  explodes slower than $t^{-(\#f+d+4)/2}$ as
  $t$ tends to 0.

\vspace{.5cm}

We can now prove the propositions stated in Section
\ref{sec:firstseries}. Proposition \ref{cor:principal1} is immediate from
Theorem \ref{thm:principal}. In the special case when $S$ is a measurable and polynomially growing function, we get
Proposition \ref{cor:f_mes_pol}:

\begin{proof}[Proof of Proposition \ref{cor:f_mes_pol}]
Multiplying (\ref{eq:principal}) by $f(y)$ and integrating in $y$ leads
to (\ref{eq:dl_croiss_pol}) with the remainder $r_n(t,x)=n^{-2}
\int_{\R^d}f(y)\pi_n(t,x,y)\: dy$. Since $\val{f(y)}\le
c'(1+\norm{y}^q)$ and $(\pi_n,n\ge 1)$ is bounded in $\cg_4(\R^{d})$, we
can find $c_1\ge 0$ and $c_2>0$ such that for all $n\ge 1$, $t\in(0,1]$
and $x\in\R^d$, $\val{r_n(t,x)}\le c_1n^{-2}t^{-(d+4)/2}
\int_{\R^d}(1+\norm{y}^q)\exp(-c_2\norm{x-y}^2/t)\: dy$. Setting
$\zeta=(y-x)/\sqrt{t}$ leads to $\val{r_n(t,x)}\le c_1n^{-2}t^{-2} \int_{\R^d}(1+\norm{x+\zeta\sqrt{t}}^{q})\exp(-c_2\norm{\zeta}^2)\: d\zeta$. To complete the proof, it remains to observe that there exists $c\ge 0$ such that for all $t\in(0,1]$ and $x,\zeta\in\R^d$, $\norm{x+\zeta\sqrt{t}}^{q}\le c(\norm x ^{q} + \norm \zeta ^{q})$.
\end{proof}

It is easy to adapt the preceding proof to get Proposition \ref{cor:exp}. In the same way, differentiating (\ref{eq:dl2}) $\alpha$
times in $x$, multiplying by $f(y)$ and integrating in $y$ leads to
Proposition \ref{cor:d/dx_exp}.

\vspace{.5cm}

Expansion (\ref{eq:principal}) should be seen as an improvement of
(\ref{eq:dl_bally-talay2}): it
allows for infinite differentiation in $x$ and $y$ and also precises the
way the coefficients explode when $t$ tends to 0. We have an analogous improvement for
expansion (\ref{eq:dl_konakov-mammen}):

\begin{theorem}\label{thm:principal_km}
Under (B) and (C), for each $i\ge 1$, there exists a bounded family $(\pi_{n,i},n\ge
1)$ in $\cg_{2i-2}(\R^d)$ and two bounded families $(\pi'_{n,i},n\ge 1)$ and $(\pi''_{n,i},n\ge 1)$ in $\cg_{2i}(\R^d)$ such that for all $j,n\ge 1$,
\begin{equation}\label{eq:dev_densite_tout_ordre}
p_n-p = \sum_{i=1}^{j-1}\frac{\pi_{n,i}}{n^i} + \sum_{i=2}^{j}\left(t-\ntn\right)^i\pi'_{n,i} + \frac{\pi''_{n,j}}{n^j}.
\end{equation}
\end{theorem}

Here and in all the sequel we use the convention that a sum over an empty
set is zero, and $\lfloor nt\rfloor$ denotes the greatest integer less
than or equal to $nt$. Expressions
involving $\lfloor nt\rfloor$ do not appear in (\ref{eq:principal}) since they
are hidden in the remainder. When $t=1$ and no differentiation is applied neither in
$x$ nor in $y$, (\ref{eq:dev_densite_tout_ordre}) boils down to the result of V. \textsc{Konakov} and
E. \textsc{Mammen} \cite{konakov-mammen}. Again note that (\ref{eq:dev_densite_tout_ordre}) is much richer in the sense that it
allows for infinite differentiation in space and also precises the
dependence on $t$. Theorem
\ref{thm:principal_km} will also be proved in Section
\ref{sec:proof_thm_principal}.

\subsection{A preliminary result}\label{sec:preliminary}

As explained in section \ref{sec:purpose}, in order to prove point (iii)
in Theorem \ref{thm:principal}, we first seek an expansion for the error operator
\begin{equation*}\label{dl_f} 
\Delta_t^n = P^n_t - P_t
\end{equation*}
where, for $f\in\cpol$ and $x\in\R^d$, we have set $P_tf(x)=\E[f(X^x_t)]$ and $P^n_tf(x)=\E\left[f(X^{n,x}_t)\right]$.
Precisely, we look for operators $C_t$ and $R^n_t$ such that $R^n_t=O(1/n^2)$ and $\Delta^n_{t}=C_t/n + R^n_t$.
The following theorem, interesting in itself, is proved in Section
\ref{sec:esp}. It can be seen as an improvement of
\cite{talay-tubaro}. It not only gives explicit formulas for $C_tf(x)$
and $R^n_tf(x)$ but also provides useful information about their
dependencies on $n,t,f$ and $x$. Note that it does not require neither
(B) nor (B') nor (C). In order to state it shortly, let us
\begin{itemize}
\item denote by $\clpol$ the space of endomorphisms of $\cpol$,
\item say that a subset $\cb\subset\cpol$ is bounded if, in (\ref{eq:C_infty_pol_Rd}), $c$ and $q$ can be chosen independently on $f\in\cb$,
\item say that $T\in\clpol$ is bounded if for all bounded $\cb\subset\cpol$, $\{Tf|f\in\cb\}$ is a bounded subset of $\cpol$,
\item denote by $\clbpol$ the space of bounded endomorphisms of $\cpol$,
\item say that a $\clbpol$-valued family $(T_i,i\in I)$ is bounded if for all bounded $\cb\subset\cpol$, $\{T_if|f\in\cb,i\in I\}$ is a bounded subset of $\cpol$,
\item say that $(T_i,i\in I)$ is a $O(h(i))$ family in $\clbpol$ if the family $(h(i)^{-1}T_i,i\in I)$ is bounded.
\end{itemize}
It is already known that, under (A), $(P_t,t\in[0,1])$ is a bounded family in $\clbpol$. A proof can de found in \cite{kusuoka-stroock}, Lemma 3.9, 
page 15. Using Lemma \ref{lem:moments_euler_A}, this proof straightforwardly adapts uniformly in $n$ so that $(P^n_t,t\in[0,1],n\ge 1)$ is also bounded in $\clbpol$. We are now in the position to state the main result of the first step:

\begin{theorem}\label{thm:principal_section_esp}
Under (A), $(\Delta^n_{t},t\in[0,1],n\ge 1)$ is a $O(t/n)$ family in
$\clbpol$, and there exists a $O(t)$ process $(C_t,t\in [0,1])$ and a
$O(1/n^2)$ family $(R^n_t,t\in [0,1],n\ge 1)$ in $\clbpol$ such that
\begin{equation*}
\Delta_{t}^n=C_t/n + R^n_t. 
\end{equation*}
Moreover, $C_t$ is explicitely given in terms of $(P_t,t\in[0,1])$ and of $L_2^*$ (see (\ref{eq:L2star})) by
\begin{equation}\label{eq:CtL2star}
C_t = \inv 2 \int_0^t P_sL_2^*P_{t-s}\: ds.
\end{equation}
\end{theorem}

Note that this theorem covers the result of D. \textsc{Talay} and L. \textsc{Tubaro}
\cite{talay-tubaro} since it implies that for any $f\in\cpol$ we can find
a $q\in\N$ such that for all
$x\in\R^d$, $t\in[0,1]$ and $n\ge 1$,
\begin{equation*}
\Delta_t^nf(x) = C_tf(x)/n + O\left(\frac{1+\norm{x}^q}{n^2} \right).
\end{equation*}
It even improves it a bit since we see that this holds under nothing more
than condition (A), whereas D. \textsc{Talay} and L. \textsc{Tubaro}
state their result under the stronger condition (B'). Note also that if we restrict ourselves to times $t$ belonging to the
discretization grid $\frak{S}_n$, we get a better control, of order $O(t/n^2)$, of the
remainder, see Remark \ref{rem:R'}.

Instead of Theorem \ref{thm:principal_section_esp}, in order to derive
Theorem \ref{thm:principal_km}, we shall need

\begin{theorem}\label{thm:principal_section_esp_km}
Under (A), there exists a sequence of differential operators
$(L^*_j,j\ge 2)$, recursively defined by (\ref{eq:recursionLxj})-(\ref{eq:defL*j}), and for each $i\ge 1$ a
$O(t/n^{i})$ family $R^i=(R^{n,i}_{t},t\in[0,1],n\ge
1)$ in $\clbpol$ such that for all $j\ge 1$,
\begin{equation}\label{eq:principal_section_esp_km}
\Delta_{t}^n = \sum_{i=2}^{j}\inv{i!n^i}\sum_{k=0}^{\nt-1}
P^n_{t^{n}_{k}}L^*_iP_{t-t^{n}_{k}} + R^{n,j}_{t} + \sum_{i=2}^{j}\frac{\left(t-\ntn\right)^i}{i!}P^n_{\ntn}L^*_iP_{t-\ntn}.
\end{equation}
\end{theorem}

Observe that the main term in (\ref{eq:principal_section_esp_km}) is
\begin{equation*}
\inv{n}\left(\inv{2n}\sum_{k=0}^{\nt-1}
P^n_{t^{n}_{k}}L^*_2P_{t-t^{n}_{k}}\right) \approx \frac{C_t}{n},
\end{equation*}
and the remainder is of order $1/n^2$. Note also that if we restrict
ourselves to times belonging to the
discretization grid $\frak{S}_n$, we get the following expansion in
$\clbpol$:
\begin{equation*}
\Delta_{\ntn}^n = \sum_{i=2}^{j}\inv{i!n^i}\sum_{k=0}^{\nt-1}
P^n_{t^{n}_{k}}L^*_iP_{\ntn-t^{n}_{k}} + O\left(\frac{t}{n^j}\right).
\end{equation*}
Theorem \ref{thm:principal_section_esp_km} is also proved in Section \ref{sec:esp}.

\subsection{Organization of the paper}

Section \ref{sec:esp} deals with the expansion for the expectation: it is dedicated to the proofs of Theorems
\ref{thm:principal_section_esp} and
\ref{thm:principal_section_esp_km}.

Section \ref{sec:density} is our second and final step. It is devoted to the proofs of Theorems
\ref{thm:principal} and \ref{thm:principal_km}. It begins with the proof
of Proposition \ref{prop:convol_esp} concerning the space convolution $*_{g,\gamma}$ in $\cg(\R^d)$.

Eventually, Section \ref{sec:appendix} is an appendix where we have
gathered useful results on the Euler scheme and technical lemmas that
are used in Sections \ref{sec:esp} and \ref{sec:density}.

\section{First step: expansion for $\E\left[f\left(X^{n,x}_t\right)\right]$}\label{sec:esp}

In this section we seek to expand $\Delta_t^nf(x)=\E\left[f\left(X^{n,x}_t\right)\right] -
\E\left[f\left(X^x_t\right)\right]$ in powers of the time step $1/n$
when $f$ is a regular function, say $f\in\cpol$. The idea is the
following. Recall the discussion preceding Theorem
\ref{thm:principal_section_esp}: under (A), both $P_t$ and $P^n_t$ are
endomorphisms of $\cpol$. In $\clpol$ we then write
\begin{equation}\label{eq:delta_nt_base}
\Delta_{t}^n = P^n_{t}-P_{t} = \sum_{k=0}^{\nt-1} P^n_{t^n_k}\Delta^n_{1/n}P_{t-t^n_{k+1}}+P^{n}_{\ntn}\Delta^n_{t-\ntn}.
\end{equation}
There is a subtle point here: $(X^{n,x}_t,t\in[0,1])$ is not a Markov
process, since the future of $X^{n,x}_t$ depends on the past value
$X^{n,x}_{\ntn}$, see (\ref{eq:def_euler}). Nevertheless, it is easy to
check by conditioning on $\cf_{t^n_k}$ that we have
$P^n_{t^n_k}P^n_s=P^n_{t^n_k+s}$ for all $s\ge 0$ - but beware: this is
different from $P^n_sP^n_{t^n_k}$ as soon as $ns$ is not an integer.

Equation (\ref{eq:delta_nt_base}) leads us to expand $\Delta_{t}^n$ for
small $t$, namely for $t\le 1/n$. This naturally involves a series of differential operators
as we shall now see.

\subsection{Operators associated with the Euler scheme}

Let us denote by $L$ the infinitesimal generator of the diffusion $X$ and by $(L^x,x\in\R^d)$ its tangent infinitesimal generator, i.e.
\begin{equation*}
L = \sum_{i=1}^{d}b_{i}\partial^{e_i}+ \frac{1}{2}\sum_{i,j=1}^{d} a_{i,j}\partial^{e_i+e_j}\qquad \mbox{and} \qquad 
L^x = \sum_{i=1}^{d}b_{i}(x)\partial^{e_i}+ \frac{1}{2}\sum_{i,j=1}^{d} a_{i,j}(x)\partial^{e_i+e_j}.
\end{equation*}
We use the convention that $L$ and
$L^x$ act on the variable $y$, so that, for instance, $L\psi(t,x,y)$ and
$L^x\psi(t,x,y)$ respectively stand for
$L\left(\psi(t,x,\cdot)\right)(y)$ and
$L^x\left(\psi(t,x,\cdot)\right)(y)$. $L^x$ is the infinitesimal
generator of the Euler scheme $(X^{n,x}_t,t\in[0,1/n])$ starting from
$x$, over the first discretization time interval: $L^x$ is built from
$L$ in the same way as $X^{n}$ is built from $X$, by freezing the drift
$b$ and the volatility $\sigma$ to their initial value on discretization
intervals. Besides, for each $x\in\R^d$ we define a sequence of differential operators $(L^x_j,j\in\N)$
by putting $L^x_0=I$ (the identity operator) and
\begin{equation}\label{eq:recursionLxj}
L^x_{j+1}=L^xL^x_j-L^x_jL,
\end{equation}
and we set
\begin{equation}\label{eq:defL*j}
L^*_jf(x)=L^x_jf(x).
\end{equation}
Observe that $L^*_1=0$. Besides, $L_2^*$ is given by (\ref{eq:L2star}) so that, under (A), $L_2^*\in\clbpol$ and there exists a family $(g^*_{2,\alpha},1\le|\alpha|\le 3)$ in $\cpol$ such that
\begin{equation}\label{eq:l*2}
L^*_2 = \sum_{1\le|\alpha|\le 3} g^*_{2,\alpha}\partial^{\alpha}.
\end{equation}
$L_2^*$ gives the exact principal error term in the expansion of
$\Delta_t^n$, see (\ref{eq:CtL2star}) and (\ref{eq:pi}). $L^*_j$ is the differential operator
appearing in (\ref{eq:principal_section_esp_km}). It does not give the exact expansion in
powers of $1/n$ but an approximated version, in the spirit of \cite{konakov-mammen}, since in (\ref{eq:principal_section_esp_km}) the
coefficients depend on $n$ - but should themselves be expanded in powers
of $1/n$. See \cite{jacod-kurtz-meleard-protter}, equations (6.35) and
(6.36), for an expression of the operators involved in the exact expansion.

Under (A), $L$ and $L^x$ belong to $\clbpol$ for each $x\in\R^d$, and, by induction, so does $L^x_j$. We can describe $L^x_j$ more precisely. Indeed, defining the powers of an operator $A$ by $A^0=I$ and $A^{j+1}=AA^j$, inductions on $j$ lead
to $L^x_j=\sum_{i=0}^j (-1)^i \binom{j}{i}\left(L^x\right)^{j-i}L^i$
and to the existence of a family $(g_{j,\alpha},h_{j,\alpha},j\in\N^*,1\le|\alpha|\le 2j)$ in $\cpol$ such that 
\begin{equation*}
 \forall x\in\R^d,\qquad \left(L^x\right)^j = \sum_{1\le|\alpha|\le 2j} g_{j,\alpha}(x)\partial^{\alpha}\qquad \mbox{and} \qquad L^j = \sum_{1\le|\alpha|\le 2j} h_{j,\alpha}\partial^{\alpha}.
\end{equation*}
Hence, for each $j\in\N^*$ one can find a family $(m_{j,\alpha},1\le |\alpha|\le 2j)$ of integers and a family $\left(g_{j,\alpha,l},h_{j,\alpha,l},1\le |\alpha|\le 2j, 1\le l\le m_{j,\alpha}\right)$ in $\cpol$ such that for all $x\in\R^d$,
\begin{equation}\label{eq:lxj}
L^x_j=\sum_{1\le |\alpha|\le 2j} 
\left(\sum_{l=1}^{m_{j,\alpha}} g_{j,\alpha,l}(x)h_{j,\alpha,l}\right)\partial^{\alpha}. 
\end{equation}

\begin{remark}\label{rem:cbo}
Note that when (B) holds, the functions $g_{j,\alpha,l}$,
$h_{j,\alpha,l}$ and $g^*_{2,\alpha}$ all belong to $\cbo$ (in fact they are polynomial in $b$, $\sigma$ and their derivatives). 
\end{remark}

We are now in the position to define a family of operators
$\Phi^j=(\Phi^{n,j}_{s,t},n\ge 1,0\le s\le t\le 1/n)$ as follows:
\begin{equation}\label{eq:def_Phi} 
\forall f\in\cpol,\quad \Phi^{n,j}_{s,t}f(x)=\E\left[L^{x}_jP_{t-s}f\left(X^{n,x}_s\right)\right].
\end{equation}
Observe that $\Phi^{n,j}_{0,t}=L_j^*P_{t}$ and that, from (\ref{eq:lxj}), 
\begin{equation}\label{eq:def_Phi_prime} 
\Phi^{n,j}_{s,t}=\sum_{1\le |\alpha|\le 2j} \sum_{l=1}^{m_{j,\alpha}} g_{j,\alpha,l}P^n_{s}(h_{j,\alpha,l}\partial^{\alpha}P_{t-s}).
\end{equation}
Boundedness is a key property of this family:

\begin{proposition}\label{prop:bornitudephipsi}
Under (A), $\Phi^j$ is a bounded family in $\clbpol$.
\end{proposition}

\begin{proof}
$(P_t,t\in[0,1])$ and $(P^n_t,t\in[0,1],n\ge 1)$ are bounded families in
$\clbpol$, see the discussion preceding Theorem \ref{thm:principal_section_esp}. Besides, multiplication by a function in $\cpol$ and differentiation are bounded operators on $\cpol$. As a sum of compositions of bounded families in $\clbpol$, $\Phi^j$ is a bounded family in $\clbpol$.
\end{proof}

The family $\Phi^j$ naturally appears when we recusively
use Itô's formula to expand $\Delta_t^n$ for small $t$, as we now explain.

\subsection{Itô expansions}

We recall (see \cite{kusuoka-stroock}, theorem 3.11, page 16) that for $f\in\cpol$, $(s,y)\mapsto
P_{t-s}f(y)$ is infinitely differentiable on $[0,t]\times\R^d$ and
\begin{equation}\label{eq:edpu}
\forall (s,y)\in[0,t]\times\R^d,\qquad (\partial_s + L)P_{t-s}f(y)=0.
\end{equation}
Since $\partial_s$ and
$L^x_j$ commute, (\ref{eq:edpu}) and the definition of $L^x_j$ imply
\begin{equation}\label{eq:lambdaj->j+1}
(\partial_s  + L^{x})L^{x}_jP_{t-s} = (L^xL^x_j-L^x_jL)P_{t-s} = L^{x}_{j+1}P_{t-s}.
\end{equation}

For a measurable family $(A_s)$ in $\clbpol$, we denote by
$\int_{t_1}^{t_2}A_s\: ds$ the element of $\clpol$ which maps $f$ to $x\mapsto\int_{t_1}^{t_2}A_sf(x)\: ds$.
The following lemma states that
$\Phi^{n,j+1}_{\cdot,t}$ is the derivative of $\Phi^{n,j}_{\cdot,t}$ on the interval $[0,t]$.

\begin{lemma}
Under (A), for all $j\in\N$, $n\ge 1$ and $0\le s\le t\le 1/n$,
\begin{equation}\label{eq:derivPhi}
\Phi^{n,j}_{s,t} = L_j^*P_{t} +
\int_{0}^{s} \Phi^{n,j+1}_{s',t}\: ds'.
\end{equation}
\end{lemma}

\begin{proof}
For $f\in\cpol$, $(s,y)\mapsto L^{x}_jP_{t-s}f(y)$ is infinitely differentiable on $[0,t]\times\R^d$ so that
we can apply Itô's formula to it and to the semimartingale $X^{n,x}$ between $0$
and $s$. Using (\ref{eq:lambdaj->j+1}) for the second equality, we get
\begin{multline*}
L^x_jP_{t-s}f\left(X^{n,x}_{s}\right)
- L^x_jP_{t}f\left(x\right) - M_s\\
= \int_{0}^{s} \left(\frac{\partial}{\partial s} +
  L^x
\right)L^x_jP_{t-s'}f\left(X^{n,x}_{s'}\right)\:
ds'
= \int_{t^{n}_{k}}^{s} L^x_{j+1}P_{t-s'}f\left(X^{n,x}_{s'}\right)\:
ds'
\end{multline*}
where $M_s=\sum_{i=1}^{d}\sum_{j=1}^{r}\sigma_{i,j}(x)
\int_{0}^{s}\partial^{e_{i}}\left(L_j^{x}P_{t-s'}f\left(X^{n,x}_{s'}\right)\right)\:
dB^{j}_{s'}$. Since $\{L_j^xP_{t-s'}f|s'\in[0,t]\}$ is bounded in $\cpol$, (\ref{eq:moments_euler}) implies that $(M_s,s\in[0,t])$ is a square-integrable martingale and thus has zero mean. Hence, taking expectations and
using (\ref{eq:def_Phi}) and Fubini's theorem, we have
\begin{equation*}
\Phi^{n,j}_{s,t}f(x)
- L_j^*P_{t}f(x)
= \int_{0}^{s} \E\left[L^x_{j+1}P_{t-s'}f\left(X^{n,x}_{s'}\right) \right] \:
ds'
= \int_{0}^{s} \Phi^{n,j+1}_{s',t}f(x)  \: ds',
\end{equation*}
which concludes the proof.
\end{proof}

For $t\in[0,1/n]$, since $\Delta^n_t=\Phi^{n,0}_{t,t}-\Phi^{n,0}_{0,t}$, by iterating (\ref{eq:derivPhi}) we get
\begin{equation}\label{eq:deltaQnkt}
\Delta^n_t = \sum_{i=2}^{j}
\frac{t^i}{i!}L_i^*P_{t} + I^{n,j+1}_{t},
\end{equation}
where
\begin{equation}\label{eq:R_nkjt}
I^{n,j+1}_{t} = \int_{0}^{t}\int_{0}^{s_{1}}\cdots\int_{0}^{s_{j}}
\Phi^{n,j+1}_{s_{j+1},t}\: ds_{j+1}\cdots ds_{2} ds_{1}.
\end{equation}
The crucial point here is that, by construction,
$L^*_1=0$ so that
the sum in (\ref{eq:deltaQnkt}) begins with $i=2$.

Injecting this in (\ref{eq:delta_nt_base}), we eventually get for all
$t\in[0,1]$ and $n\ge 1$
\begin{equation}\label{eq:delta_nt_dev}
\Delta_{t}^n = \sum_{i=2}^{j}\inv{i!n^i}\sum_{k=0}^{\nt-1}
P^n_{t^{n}_{k}}L^*_iP_{t-t^{n}_{k}} + R^{n,j}_{t} +
\sum_{i=2}^{j}\frac{\left(t-\ntn\right)^i}{i!}P^n_{\ntn}L^*_iP_{t-\ntn},
\end{equation}
where
\begin{equation}\label{eq:R_njt}
R^{n,j}_{t} =
\sum_{k=0}^{\nt-1}P^n_{t^n_k}I^{n,j+1}_{1/n}P_{t-t^n_{k+1}} +
P^n_{\ntn}I^{n,j+1}_{t-\ntn}.
\end{equation}

From Proposition \ref{prop:bornitudephipsi},
$(I^{n,j+1}_{t},n\ge 1,t\in[0,1/n])$ is a
$O(t^{j+1})$ family in $\clbpol$. Recalling the boundedness of
$(P_t,t\in[0,1])$ and $(P^n_t,t\in[0,1],n\ge 1)$, we get that the family $R^i=(R^{n,i}_{t},t\in[0,1],n\ge
1)$ is $O(t/n^{i})$ in $\clbpol$. Theorem
\ref{thm:principal_section_esp_km} is thus proved. We are now in good position to prove Theorem \ref{thm:principal_section_esp}.

\subsection{Proof of Theorem \ref{thm:principal_section_esp}}

In the particular case when $j=1$, (\ref{eq:delta_nt_dev}) reads
$\Delta_{t}^n = R^{n,1}_t$ so that we have proved that $(\Delta_{t}^n,t\in[0,1],n\ge 1)$ is $O(t/n)$
in $\clbpol$, which was the first statement of Theorem \ref{thm:principal_section_esp}.

In the particular case when $j=2$, if we set
\begin{eqnarray}\label{eq:CtRnt}
C_t &=& \frac{1}{2} \int_{0}^{t}P_{s}L^*_2P_{t-s}\: ds, \label{eq:Ct}\\
A^{n}_{1,t} &=& \frac{1}{2n}\left(\inv n \sum_{k=0}^{\nt-1}
P_{t^{n}_{k}}L^*_2P_{t-t^{n}_{k}} -
\int_{0}^{t}P_{s}L^*_2P_{t-s}\: ds \right), \label{eq:R^n_1,t} \\
A^{n}_{2,t} &=& \frac{1}{2n^2}\sum_{k=0}^{\nt-1}
\left(P^n_{t^{n}_{k}}L^*_2P_{t-t^{n}_{k}} -
P_{t^{n}_{k}}L^*_2P_{t-t^{n}_{k}} \right), \label{eq:R^n_2,t} \\
A^{n}_{3,t} &=& R^{n,2}_{t} + \frac{\left(t-\ntn\right)^2}{2}P^n_{\ntn}L^*_2P_{t-\ntn} \label{eq:An3t}\\
R^n_t &=& A^{n}_{1,t} + A^{n}_{2,t} + A^{n}_{3,t}, \label{eq:Rnt}
\end{eqnarray}
equation (\ref{eq:delta_nt_dev}) reads
\begin{equation}\label{eq:delta_nt_dev_2}
\Delta_{t}^n = C_t/n + R^n_t. 
\end{equation}
As a composition of bounded families, $(P_{s}L^*_2P_{t-s},0\le s\le t\le
1)$ is a bounded family in $\clbpol$, so that $(C_t,t\in [0,1])$ is
$O(t)$ in $\clbpol$. It remains to prove that $(R^n_t,t\in [0,1],n\ge
1)$ is $O(1/n^2)$ in $\clbpol$. We have already proved that it is true of
$(R^{n,2}_{t},t\in [0,1],n\ge
1)$. It is obviously also true of $(\left(t-\ntn\right)^2P^n_{\ntn}L^*_2P_{t-\ntn},t\in [0,1],n\ge
1)$, so that $(A^n_{3,t},t\in [0,1],n\ge
1)$ is $O(1/n^2)$ in $\clbpol$.

For $(A^n_{1,t},t\in [0,1],n\ge
1)$, observe that, if we
set $L^{\#}_3=LL^*_2-L^*_2L\in\clbpol$, as $\partial_sP_s=LP_s=P_sL$, we have $\partial_sP_{s}L^*_2P_{t-s}=P_{s}LL^*_2P_{t-s}-P_{s}L^*_2LP_{t-s}=P_{s}L^{\#}_3P_{t-s}$. Hence the family $(P_{t^{n}_{k}}L^*_2P_{t-t^{n}_{k}} -
P_{s}L^*_2P_{t-s},t\in[0,1],n\ge 1,k\in\{0,\ldots,\nt-1\},s\in[t^{n}_{k},t^{n}_{k+1}])$ satisfies
\begin{equation}\label{eq:Rn1t1}
P_{t^{n}_{k}}L^*_2P_{t-t^{n}_{k}} - P_{s}L^*_2P_{t-s}=-\int_{t^n_k}^s P_{u}L^{\#}_3P_{t-u}\: du
\end{equation}
and thus is $O(1/n)$ in $\clbpol$. As a consequence,
\begin{equation}\label{eq:Rn1t2}
A^{n}_{1,t} = \frac{1}{2n}\sum_{k=0}^{\nt-1}
\int_{t^n_k}^{t^n_{k+1}}\left(P_{t^{n}_{k}}L^*_2P_{t-t^{n}_{k}} -
P_{s}L^*_2P_{t-s}\right)\: ds - \frac{1}{2n}\int_{\ntn}^{t}P_{s}L^*_2P_{t-s}\: ds
\end{equation}
is $O(1/n^2)$ in $\clbpol$.

As for $(A^n_{2,t},t\in [0,1],n\ge
1)$, note that $P^n_{t^{n}_{k}}L^*_2P_{t-t^{n}_{k}}-P_{t^{n}_{k}}L^*_2P_{t-t^{n}_{k}}=\Delta^n_{t^n_k}L^*_2P_{t-t^n_k}$.
Since $(\Delta_{t}^n,t\in[0,1],n\ge 1)$ is $O(1/n)$
in $\clbpol$, so is the family $(P^n_{t^{n}_{k}}L^*_2P_{t-t^{n}_{k}}-P_{t^{n}_{k}}L^*_2P_{t-t^{n}_{k}},t\in[0,1],n\ge 1,k\in\{0,\ldots,\nt-1\})$, as the composition of a bounded family by a
$O(1/n)$ family in $\clbpol$. This completes the proof of Theorem \ref{thm:principal_section_esp}.

\begin{remark}\label{rem:R'}
It is noteworthy that the family $(R'^n_t,t\in [0,1],n\ge 1)$ defined by
\begin{equation*}
R'^n_t=R^n_t+ \inv{2n}\int_{\ntn}^{t}P_{s}L^*_2P_{t-s}\: ds - \frac{\left(t-\ntn\right)^2}{2}P^n_{\ntn}L^*_2P_{t-\ntn}
\end{equation*}
is $O(t/n^2)$ in $\clbpol$. In particular, $(R^n_{\ntn},t\in [0,1],n\ge 1)$ is $O(t/n^2)$ in $\clbpol$.
\end{remark}


\section{Second step: expansion for the density of $X^{n,x}_t$}\label{sec:density}

This section is devoted to the proofs of Theorems \ref{thm:principal} and
\ref{thm:principal_km}.

\subsection{Space convolutions}

We begin by proving Proposition
\ref{prop:convol_esp} which is the key argument. Recall the definitions
of Section \ref{sec:def_spaces}. Let $\cb_1$ and $\cb_2$
be two bounded subsets of $\cg(\R^{d})$, $g\in\cbo$ and
$\gamma\in\N^d$. We want to prove that
\begin{itemize}
\item[(i)] $\{\pi_1*_{g,\gamma}\pi_2|\pi_1\in\cb_1,\pi_2\in\cb_2\}$ is a
  bounded subset of $\ch_{\val{\gamma}}(\R^{d})$,
\item[(ii)] $\{\pi_1*_{g,0}\pi_2|\pi_1\in\cb_1,\pi_2\in\cb_2\}$ is a bounded subset of $\ch(\R^{d})$.
\end{itemize}

The functions $\pi_1*_{g,\gamma}\pi_2$ depend on
$(s,t,x,y)$. We shall proceed differently depending on $s$ is small or
large with respect to $t$. The main trick is to integrate by parts in
the latter case, so that the derivatives should always rest on the regularizing
part of the integral. This is analogous to V. \textsc{Bally} and
D. \textsc{Talay}'s use of Malliavin calculus integration by parts
formula \cite{bally-talay1}. This is the reason why we partition the unit
triangle $\ct_1$ into $\ct_1^-=\{(s,t)\in\ct_1|0<s\le t/2\}$ and
$\ct_1^+=\{(s,t)\in\ct_1|t/2<s<t\}$, and, for $\epsilon=\pm$, we define
$\left(\pi_1*_{g,\gamma}\pi_2\right)_{\epsilon}(s,t,x,y)=\ind_{\ct_1^{\epsilon}}(s,t)\left(\pi_1*_{g,\gamma}\pi_2\right)(s,t,x,y)$.
We then have $\pi_1*_{g,\gamma}\pi_2 = (\pi_1*_{g,\gamma}\pi_2)_- +(\pi_1*_{g,\gamma}\pi_2)_+$.

Before proving Proposition \ref{prop:convol_esp} and for the sake of clarity, let us state apart the following technical lemma, whose proof is a straightforward application of Lebesgue's dominated convergence theorem:
\begin{lemma}\label{lem:convol_esp}
Let $l\in\Z$, $(\chi_i,i\in I)$ be a family of measurable functions mapping $\ct_1\times\R^{d}\times\R^{d}\times\R^{d}$ into $\R$ such that
\begin{itemize}
\item for all $i\in I$, $(s,t)\in\ct_1$ and $\zeta\in\R^d$, $\chi_i(s,t,\cdot,\cdot,\zeta)$ is infinitely differentiable,
\item for all $\alpha,\beta\in \N^{d}$, there exists two constants
  $c_1\ge 0$ and $c_2>0$ such that for all $i\in I$, $(s,t)\in\ct_1$
  and $x,y,\zeta\in\R^d$,
\begin{equation}\label{eq:def_cil1}
\val{\partial^{\alpha}_x\partial^{\beta}_y\chi_i(s,t,x,y,\zeta)} \le c_1 t^{-\left(\val\alpha+\val\beta+d+l\right)/2}\exp\left(-c_2\norm{x-y}^2/t-c_2\norm{\zeta}^2\right),
\end{equation}
\end{itemize}
and let us define $\ci(\chi_i)(s,t,x,y)=\int_{\R^d}\chi_i(s,t,x,y,\zeta)\: d\zeta$. Then $\{\ci(\chi_i)|i\in I\}$ is a bounded subset of $\ch_l(\R^d)$.
\end{lemma}

\begin{proof}[Proof of Proposition \ref{prop:convol_esp}-(i).]
It is enough to show that both $\cb_{\epsilon}\equiv\{(\pi_1*_{g,\gamma}\pi_2)_{\epsilon}|\pi_1\in\cb_1,\pi_2\in\cb_2\}$ are bounded.

{\bf Step 1.} Let us first treat $\cb_-$, i.e. the case when $s$ is small. After the change of variables $z=x+\zeta\sqrt{s}$, we get $(\pi_1*_{g,\gamma}\pi_2)_-=\ci(\chi^-_{\pi_1,\pi_2})$ with
\begin{equation*}
 \chi^-_{\pi_1,\pi_2}(s,t,x,y,\zeta) = \ind_{\ct_1^-}(s,t)s^{d/2}g(x+\zeta\sqrt{s})\pi_1(s,x,x+\zeta\sqrt{s})\partial^{\gamma}_2\pi_2(t-s,x+\zeta\sqrt{s},y).
\end{equation*}
It is enough to check that the family $\left(\chi^-_{\pi_1,\pi_2},(\pi_1,\pi_2)\in\cb_1\times\cb_2\right)$ satisfies the assumptions of Lemma \ref{lem:convol_esp} with $l=\val\gamma$. The first point is obvious. In order to check the second one, let us fix $\alpha,\beta\in \N^{d}$. According to Leibniz's formula, $\partial^{\alpha}_x\partial^{\beta}_y\chi_{\pi_1,\pi_2}(s,t,x,y,\zeta)$ can be written as a weighted sum of terms of the form
\begin{multline*}
 \chi^{-,\alpha_1,\alpha_2,\alpha_3}_{\pi_1,\pi_2}(s,t,x,y,\zeta) = \ind_{\ct_1^-}(s,t)s^{d/2}\partial^{\alpha_1}g(x+\zeta\sqrt{s}) \\ \partial^{\alpha_2}_x\left(\pi_1(s,x,x+\zeta\sqrt{s})\right)\partial^{\gamma+\alpha_3}_2\partial^{\beta}_3\pi_2(t-s,x+\zeta\sqrt{s},y),
\end{multline*}
with $\val{\alpha_1}+\val{\alpha_2}+\val{\alpha_3}=\val\alpha$, so that
in order to check (\ref{eq:def_cil1}) it is enough to show that for each
such $(\alpha_1,\alpha_2,\alpha_3)$ one can find $c_1\ge 0$ and $c_2>0$
such that for all $(\pi_1,\pi_2)\in\cb_1\times\cb_2$, $(s,t)\in\ct_1$
and $x,y,\zeta\in\R^d$,
$|\chi^{-,\alpha_1,\alpha_2,\alpha_3}_{\pi_1,\pi_2}(s,t,x,y,\zeta)|$ is
less than the r.h.s. of (\ref{eq:def_cil1}), with $l=\val\gamma$. Now,
$\cb_1$ and $\cb_2$ are bounded subsets of $\cg(\R^d)$ so that from
(\ref{eq:def_cg1})-(\ref{eq:def_cg2}) one can find $c_3,c_5\ge 0$ and $c_4>0$ such that for all $(\pi_1,\pi_2)\in\cb_1\times\cb_2$, $(s,t)\in\ct_1$ and $x,y,\zeta\in\R^d$,
\begin{equation*}
\val{\partial^{\alpha_2}_x\left(\pi_1(s,x,x+\zeta\sqrt{s})\right)} \le c_3s^{-d/2}\exp\left(-c_4\norm{\zeta}^2\right)
\end{equation*}
and
\begin{multline*}
\ind_{\ct_1^-}(s,t)\val{\partial^{\gamma+\alpha_3}_2\partial^{\beta}_3\pi_2(t-s,x+\zeta\sqrt{s},y)} \\ 
\le \ind_{\ct_1^-}(s,t)c_3 (t-s)^{-\left(\val{\alpha_3}+\val\beta+\val\gamma+d\right)/2}\exp\left(-c_4\norm{x-y+\zeta\sqrt{s}}^2/(t-s)\right)\\
\le \ind_{\ct_1^-}(s,t)c_5 t^{-\left(\val\alpha+\val\beta+\val\gamma+d\right)/2}\exp\left(-c_4\norm{x-y+\zeta\sqrt{s}}^2/t\right)
\end{multline*}
where, for the last inequality, we have used the fact that when $(s,t)\in\ct_1^-$,
$t/2\le t-s\le t\le 1$. Now, using the fact that $\norm{x-z}^2\ge\norm x^2/2-\norm z^2$ for all $x,z\in\R^d$, we see that for all $(s,t)\in\ct_1^-$, $\norm{\zeta}^2+\norm{x-y+\zeta\sqrt{s}}^2/t\ge(\norm{x-y}^2/t+\norm{\zeta}^2)/2$. Since $g\in\cbo$, we can eventually find $c_1\ge 0$ and $c_2>0$ such that for all $(\pi_1,\pi_2)\in\cb_1\times\cb_2$, $(s,t)\in\ct_1$ and $x,y,\zeta\in\R^d$,
\begin{equation*}
\val{\chi^{-,\alpha_1,\alpha_2,\alpha_3}_{\pi_1,\pi_2}(s,t,x,y,\zeta)} \le c_1 t^{-\left(\val\alpha+\val\beta+d+\val{\gamma}\right)/2}\exp\left(-c_2\norm{x-y}^2/t-c_2\norm{\zeta}^2\right),
\end{equation*}
which completes Step 1.

{\bf Step 2.} Let us now treat $\cb_+$, i.e. the case when $s$ is large. After $\val\gamma$ integrations by parts, we have
\begin{equation*}
(\pi_1*_{g,\gamma}\pi_2)_+(s,t,x,y)=\ind_{\ct_1^+}(s,t)\int_{\R^{d}}\partial^{\gamma}_z(g(z)\pi_1(s,x,z))\pi_2(t-s,z,y)\: dz.
\end{equation*}
Using Leibniz's formula and making the change of variables $z=y-\zeta\sqrt{t-s}$, we get that $(\pi_1*_{g,\gamma}\pi_2)_+$ is a weighted sum of terms of the form $\ci(\chi^{+,\gamma_1,\gamma_2}_{\pi_1,\pi_2})$ with
\begin{multline*}
 \chi^{+,\gamma_1,\gamma_2}_{\pi_1,\pi_2}(s,t,x,y,\zeta) = \ind_{\ct_1^+}(s,t)(t-s)^{d/2}\partial^{\gamma_1}g(y-\zeta\sqrt{t-s})\\ \partial^{\gamma_2}_3\pi_1(s,x,y-\zeta\sqrt{t-s})\pi_2(t-s,y-\zeta\sqrt{t-s},y)
\end{multline*}
and $\val{\gamma_1}+\val{\gamma_2}=\val{\gamma}$, so that we are now in the position to apply the same arguments as
in Step 1 and get that the family $(\chi^{+,\gamma_1,\gamma_2}_{\pi_1,\pi_2},(\pi_1,\pi_2)\in\cb_1\times\cb_2)$ satisfies the assumptions of Lemma \ref{lem:convol_esp} with $l=\val\gamma$, which completes the proof.
\end{proof}

\begin{proof}[Proof of Proposition \ref{prop:convol_esp}-(ii).]
From (i), we know that
$\{\pi_1*_{g,0}\pi_2|\pi_1\in\cb_1,\pi_2\in\cb_2\}$ is a bounded subset
of $\ch_0(\R^{d})$. It remains to prove that (\ref{eq:def_ch2}) holds
for $\rho=\pi_1*_{g,0}\pi_2$ with constants $c_1$ and $c_2$ which do not
depend on $(\pi_1,\pi_2)\in\cb_1\times\cb_2$. As in the proof of
Proposition \ref{prop:convol_esp}-(i), we treat $(\pi_1*_{g,0}\pi_2)_-$
and $(\pi_1*_{g,0}\pi_2)_+$ separately but analogously. That is, after
integrating by parts, the term $(\pi_1*_{g,0}\pi_2)_+$ can be treated in the
same way as $(\pi_1*_{g,0}\pi_2)_-$. Thus we shall only deal with the
latter term. We have $(\pi_1*_{g,0}\pi_2)_-=\ci(\chi^-_{\pi_1,\pi_2})$ with
\begin{equation*}
 \chi^-_{\pi_1,\pi_2}(s,t,x,y,\zeta) = \ind_{\ct_1^-}(s,t)s^{d/2}g(x+\zeta\sqrt{s})\pi_1(s,x,x+\zeta\sqrt{s})\pi_2(t-s,x+\zeta\sqrt{s},y).
\end{equation*}
Then we write $\partial^{\alpha}_x\left(\chi^-_{\pi_1,\pi_2}\left(s,t,x,x+y\sqrt{t},\zeta\right)\right)$ as a weighted sum of terms of the form
\begin{multline*}
 \tilde{\chi}^{-,\alpha_1,\alpha_2,\alpha_3}_{\pi_1,\pi_2}(s,t,x,y,\zeta) = \ind_{\ct_1^-}(s,t)s^{d/2}\partial^{\alpha_1}g(x+\zeta\sqrt{s}) \\ \partial^{\alpha_2}_x\left(\pi_1(s,x,x+\zeta\sqrt{s})\right)\partial^{\alpha_3}_x\left(\pi_2(t-s,x+\zeta\sqrt{s},x+y\sqrt{t})\right),
\end{multline*}
with $\val{\alpha_1}+\val{\alpha_2}+\val{\alpha_3}=\val\alpha$. Then we
use (\ref{eq:def_cg2}) twice and the same arguments as in the preceding
proof to get $c_1\ge 0$ and $c_2>0$ such that for all
$(\pi_1,\pi_2)\in\cb_1\times\cb_2$, $(s,t)\in\ct_1$ and
$x,y,\zeta\in\R^d$,
$|\tilde{\chi}^{-,\alpha_1,\alpha_2,\alpha_3}_{\pi_1,\pi_2}(s,t,x,y,\zeta)|\le
c_1t^{-d/2}\exp(-c_2\norm{y}^2-c_2\norm{\zeta}^2)$, and an obvious
adaptation of Lemma \ref{lem:convol_esp} completes the proof.
\end{proof}

\subsection{Proof of Theorems \ref{thm:principal} and \ref{thm:principal_km}}\label{sec:proof_thm_principal}

In this section, we assume (B) and (C). We first want to prove Theorem
\ref{thm:principal}. We recall that statement (i) is already known, see
\cite{friedman1}, theorem 7, page 260. The next lemma is statement (ii). 

\begin{lemma}\label{lem:pn_bornee}
Under (B) and (C), for all $t\in(0,1]$, $n\ge 1$ and $x\in\R^d$, $X_{t}^{n,x}$ has a density $p_n(t,x,\cdot)$ and $(p_n,n\ge 1)$ is a bounded sequence in $\cg(\R^d)$.
\end{lemma}

\begin{proof}
It is known that for all $n\ge 1$, $k\in\{1,\ldots,n\}$ and $x\in\R^d$, $X_{t^n_k}^{n,x}$ has a density $p_{n,k}(x,\cdot)$ such that $p_{n,k}$ is infinitely differentiable and satisfies (\ref{eq:def_cg1})-(\ref{eq:def_cg2}) with $t=t^n_k$ and two constants $c_1$ and $c_2$ which do not depend on $n$ and $k$ (see the proof of theorem 1.1, page 278, in \cite{konakov-mammen}). Since $\ntn\ge t/2$ for all $t\ge 1/n$, this shows that the sequence $(\tilde{p}_n,n\ge 1)$ defined by $\tilde{p}_n(t,x,y)=\ind_{\{nt\ge 1\}}p_{n,\nt}(x,y)$ is bounded in $\cg(\R^d)$.
If we denote by $\Gamma(t,x,\cdot)$ the density of $x+b(x)t+\sigma(x)B_t$ ($t\in (0,1]$), we observe that when $k\in\{1,\ldots,n-1\}$ and $t\in(t^n_k,t^n_{k+1})$, $X_{t}^{n,x}$ has the density $p_n(t,x,\cdot)=\int_{\R^d}p_{n,k}(x,z)\Gamma(t-t^n_k,z,\cdot)\: dz=(\tilde{p}_n*_{1,0}\Gamma)(t^n_k,t,x,\cdot)$. Hence, for all $t\in(0,1]$, $n\ge 1$ and $x\in\R^d$, $X_{t}^{n,x}$ has the density
\begin{equation*}
p_n(t,x,\cdot) = \left\{
          \begin{array}{ll}
            p_{n,k}(x,\cdot) & \qquad \mathrm{if}\quad t=t^n_k,\: k\in\{1,\ldots,n\}, \\
            \Gamma(t,x,\cdot) & \qquad \mathrm{if}\quad t\in(0,t^n_1), \\
            (\tilde{p}_n*_{1,0}\Gamma)(t^n_k,t,x,\cdot) & \qquad \mathrm{if}\quad t\in(t^n_k,t^n_{k+1}),\: k\in\{1,\ldots,n-1\}.
          \end{array}
        \right.
\end{equation*}
Observing that $\Gamma\in\cg(\R^d)$ and applying Proposition \ref{prop:convol_esp}-(ii), we get that $(p_n,n\ge 1)$ is a bounded sequence in $\cg(\R^{d})$.
\end{proof}

We shall now prove statement (iii) of Theorem \ref{thm:principal}. Recall
(\ref{eq:delta_nt_dev_2}). We want to make explicit $C_t$ and $R^n_t$ as integral operators on
$\R^d$. To this end, note that, applying recursively Lebesgue's dominated convergence theorem, we have that for all $t\in(0,1]$, $f\in\cpol$, $x\in\R^d$ and $\alpha\in\N^d$,
\begin{equation}\label{eq:deriv_Ptf}
\partial^{\alpha}P_tf(x) = \int_{\R^{d}} f(y)\partial^{\alpha}_{2}p(t,x,y)\: dy. 
\end{equation}

The next lemma explicits $C_t$ as an integral operator. The function
$\pi$ which appears there should be thought of as the kernel of $C$.

\begin{lemma}\label{lem:C_t}
Under (B) and (C), there exists $\pi\in\cg_1(\R^d)$, given by (\ref{eq:pi}), such that for all $t\in(0,1]$, $f\in\cpol$ and $x\in\R^d$,
\begin{equation*}
C_tf(x)=\int_{\R^d}f(y)\pi(t,x,y)\: dy.
\end{equation*}
\end{lemma}

\begin{proof}
Using (\ref{eq:Ct}) for the first equality,
(\ref{eq:l*2}) for the third one and (\ref{eq:deriv_Ptf}) for the fourth
one, we have
\begin{eqnarray*}
2C_tf(x) &=& \int_{0}^{t}P_{s}L^*_2P_{t-s}f(x)\: ds \\
         &=& \int_{0}^{t}\int_{\R^d}p(s,x,z)L^*_2P_{t-s}f(z)\: dzds \\
         &=& \sum_{1\le\val\alpha\le3}\int_{0}^{t}\int_{\R^d}g^*_{2,\alpha}(z)p(s,x,z)\partial^{\alpha}P_{t-s}f(z)\: dzds \\
         &=& \sum_{1\le\val\alpha\le3}\int_{0}^{t}\int_{\R^d}\int_{\R^d}f(y)g^*_{2,\alpha}(z)p(s,x,z)\partial^{\alpha}_2p(t-s,z,y)\: dydzds.
\end{eqnarray*}
Using Fubini's theorem, we see that to complete the proof it is enough to show that the function $\pi$ defined by
\begin{equation}\label{eq:def_pi}
\pi(t,x,y) = \inv 2 \sum_{1\le\val\alpha\le3}\int_{0}^{t} (p*_{g^*_{2,\alpha},\alpha}p)(s,t,x,y)\: ds
\end{equation}
belongs to $\cg_1(\R^d)$. Now, $p\in\cg(\R^d)$ and, from Remark \ref{rem:cbo}, $g^*_{2,\alpha}\in\cbo$ so that we can apply Proposition \ref{prop:convol_esp}-(i): $p*_{g^*_{2,\alpha},\alpha}p\in\ch_{\val\alpha}(\R^d)$. In particular, $\int_{0}^{\cdot} (p*_{g^*_{2,\alpha},\alpha}p)(s,\cdot,\cdot,\cdot)\: ds \in\cg_{\val{\alpha}-2}(\R^d)$. Since $\val{\alpha}\le 3$ and by monotonicity of $(\cg_l(\R^d),l\in\Z)$, we finally get that $\pi\in\cg_1(\R^d)$. To complete the proof, note that (\ref{eq:def_pi}) can be rewritten as (\ref{eq:pi}).
\end{proof}

We have a similar representation for $A^n_{1,t}$, recall (\ref{eq:R^n_1,t}). We say that a sequence $(\pi^n,n\ge 1)$ is $O(1/n^j)$ in $\cg_l(\R^d)$ if $(n^j\pi^n,n\ge 1)$ is bounded in $\cg_l(\R^d)$.

\begin{lemma}\label{lem:R^n_1,t}
Under (B) and (C), there exists a $O(1/n^2)$ sequence $(\pi^n_1,n\ge 1)$ in $\cg_3(\R^d)$ such that for all $t\in(0,1]$, $f\in\cpol$ and $x\in\R^d$,
\begin{equation*}
A^n_{1,t}f(x)=\int_{\R^d}f(y)\pi^n_1(t,x,y)\: dy.
\end{equation*}
\end{lemma}

\begin{proof}
Recall (\ref{eq:Rn1t1}). From Remark \ref{rem:cbo}, there is a family $(g^{\#}_{3,\alpha},1\le\val\alpha\le 4)$ in $\cbo$ such that $L^{\#}_3=\sum_{1\le\val\alpha\le 4}g^{\#}_{3,\alpha}\partial^{\alpha}$, so that, using (\ref{eq:deriv_Ptf}), we have
\begin{multline*}
(P_{t^{n}_{k}}L^*_2P_{t-t^{n}_{k}} - P_{s}L^*_2P_{t-s})f(x) \\ =
 -\sum_{1\le\val\alpha\le 4}\int_{t^n_k}^s\int_{\R^d}\int_{\R^d}f(y) g^{\#}_{3,\alpha}(z)p(u,x,z)\partial^{\alpha}_2p(t-u,z,y)\: dydzdu.
\end{multline*}
Using (\ref{eq:Rn1t2}), we get $A^n_{1,t}f(x)=\int_{\R^d}f(y)\pi^n_1(t,x,y)\: dy$ with $\pi^n_1=\pi^n_{1,1}+\pi^n_{1,2}$ and
\begin{eqnarray*}
\pi^{n}_{1,1}(t,x,y) &=& -\frac{1}{2n}\sum_{1\le\val\alpha\le 4}\sum_{k=0}^{\nt-1}
\int_{t^n_k}^{t^n_{k+1}}\int_{t^n_k}^s \left(p*_{g^{\#}_{3,\alpha},\alpha}p\right)(u,t,x,y)\: duds, \\
\pi^{n}_{1,2}(t,x,y) &=& -\frac{1}{2n}\sum_{1\le\val\alpha\le 3}\int_{\ntn}^{t} (p*_{g^*_{2,\alpha},\alpha}p)(s,t,x,y)\: ds.
\end{eqnarray*}
Now Proposition \ref{prop:convol_esp}-(i) states that
$p*_{g^{\#}_{3,\alpha},\alpha}p$ and $p*_{g^*_{2,\alpha},\alpha}p$
belong to $\ch_{\val\alpha}(\R^d)$. Hence
$(\int_{t^n_k}^{t^n_{k+1}}\int_{t^n_k}^s
(p*_{g^{\#}_{3,\alpha},\alpha}p)(u,\cdot,\cdot,\cdot)\: duds,n\ge
1,k\in\{0,\ldots,n-1\})$ is $O(1/n^2)$ in $\cg_{\val{\alpha}}(\R^d)$ and
$(\int_{\lfloor n\cdot\rfloor/n}^{\cdot}
(p*_{g^*_{2,\alpha},\alpha}p)(s,\cdot,\cdot,\cdot)\: ds,n\ge 1)$ is
$O(1/n)$ in $\cg_{\val{\alpha}}(\R^d)$. As a consequence,
$(\pi^{n}_{1,1},n\ge 1)$ is $O(1/n^2)$ in $\cg_2(\R^d)$ and
$(\pi^{n}_{1,2},n\ge 1)$ is $O(1/n^2)$ in $\cg_3(\R^d)$. Eventually, $(\pi^{n}_{1},n\ge 1)$ is $O(1/n^2)$ in $\cg_3(\R^d)$.
\end{proof}

We shall now prove an analogous lemma for $A^n_{2,t}$. 

\begin{lemma}\label{lem:R^n_2,t}
Under (B) and (C), there exists a $O(1/n^2)$ sequence $(\pi^n_2,n\ge 1)$ in $\cg_3(\R^d)$ such that for all $t\in(0,1]$, $f\in\cpol$ and $x\in\R^d$,
\begin{equation*}
A^n_{2,t}f(x)=\int_{\R^d}f(y)\pi^n_2(t,x,y)\: dy.
\end{equation*}
\end{lemma}

\begin{proof}
Since $P^n_{t^{n}_{k}}L^*_2P_{t-t^{n}_{k}} =
P_{t^{n}_{k}}L^*_2P_{t-t^{n}_{k}}$ when $k=0$, (\ref{eq:R^n_2,t}) reads 
\begin{eqnarray*}
2n^2 A^{n}_{2,t}f(x) &=& \sum_{k=1}^{\nt-1}
 (P^n_{t^n_k}-P_{t^n_k})L^*_2P_{t-t^n_k}f(x) \\
   & = & \sum_{k=1}^{\nt-1}
   \int_{\R^d}(p_n-p)(t^n_k,x,z)L^*_2P_{t-t^n_k}f(z)\: dz \\
&=&\sum_{1\le\val\alpha\le 3}\sum_{k=1}^{\nt-1} \int_{\R^d}(p_n-p)(t^n_k,x,z)g^*_{2,\alpha}(z)\partial^{\alpha}P_{t-t^n_k}f(z)\: dz \\
                & = & \sum_{1\le\val\alpha\le 3}\sum_{k=1}^{\nt-1} \int_{\R^d} \int_{\R^d}(p_n-p)(t^n_k,x,z)g^*_{2,\alpha}(z)f(y)\partial^{\alpha}_2p(t-t^n_k,z,y)\: dydz
\end{eqnarray*}
where we have used (\ref{eq:l*2}) for the third equality and
(\ref{eq:deriv_Ptf}) for the fourth one. From Remark \ref{rem:cbo}, $g^*_{2,\alpha}\in\cbo$ so that to complete the proof
it is enough to show that whenever $g\in\cbo$ and $\alpha\in\N^d$, the
sequence $(\pi^n,n\ge 1)$ defined by
\begin{equation*}
\pi^n(t,x,y) = \!\sum_{k=1}^{\nt-1} \! \int_{\R^d}
(p_n-p)(t^n_k,x,z)g(z)\partial^{\alpha}_2p(t-t^n_k,z,y)\: dz = \!\sum_{k=1}^{\nt-1}((p_n-p)*_{g,\alpha}p)(t^n_k,t,x,y)
\end{equation*}
is bounded in $\cg_{\val\alpha}(\R^d)$. And to do so, it is enough to show that the
sequence $(\rho^n_{t^n_k},n\ge 2,k\in\{1,\ldots,n-1\})$ defined by 
\begin{equation*}
\rho^n_{t^n_k}(t,x,y)=\ind_{\ct_1}(t^n_k,t)\left((p_n-p)*_{g,\alpha}p\right)(t^n_k,t,x,y)
\end{equation*}
is $O(1/n)$ in $\cg_{\val\alpha+2}(\R^d)$. Let us write $\rho^{n,-}_{t^n_k}(t,x,y) =
\ind_{\ct_1^-}(t^n_k,t)\rho^n_{t^n_k}(t,x,y)$ and $\rho^{n,+}_{t^n_k}(t,x,y) =
\ind_{\ct_1^+}(t^n_k,t)\rho^n_{t^n_k}(t,x,y)$ so that
$\rho^n_{t^n_k}=\rho^{n,-}_{t^n_k}+\rho^{n,+}_{t^n_k}$.

Let us first prove that $(\rho^{n,-}_{t^n_k},n\ge
2,k\in\{1,\ldots,n-1\})$ is $O(1/n)$ in $\cg_{\val\alpha+2}(\R^d)$. The sequence
$(\pi_{t^n_k},n\ge 2,k\in\{1,\ldots,n-1\})$ defined by $\pi_{t^n_k}(t,x,y)=\ind_{\ct_1^-}(t^n_k,t)g(x)\partial^{\alpha}_2p(t-t^n_k,x,y)$ is bounded in
$\cg_{\val\alpha}(\R^d)$, since $t-t^n_k\ge t/2$ when
$(t^n_k,t)\in\ct_1^-$. Now note that
$\rho^{n,-}_{t^n_k}=P^n_{t^n_k}\pi_{t^n_k}-P_{t^n_k}\pi_{t^n_k}=\Delta^n_{t^n_k}\pi_{t^n_k}$
(see (\ref{eq:def_P_cgl}) in the appendix for the definition of
$P^n_s\pi$, $P_s\pi$ and $\Delta^n_{s}\pi$ when $\pi\in\cg_l(\R^d)$). Thus, from (\ref{eq:delta_nt_dev})-(\ref{eq:R_njt}) and (\ref{eq:R_nkjt}) applied with $j=1$,
\begin{equation*}
\rho^{n,-}_{t^n_k}=\sum_{m=0}^{k-1}\int_{0}^{1/n}\int_{0}^{s_1}P^n_{t^n_m}\Phi^{n,2}_{s_2,1/n}P_{t^n_k-t^n_{m+1}}\pi_{t^n_k}\:
ds_2ds_1.
\end{equation*}
Proposition \ref{prop:operators_on_cgl} in the appendix states that the
family $(P^n_{t^n_m}\Phi^{n,2}_{s,1/n}P_{t^n_k-t^n_{m+1}}\pi_{t^n_k},n\ge 2,k\in\{1,\ldots,n-1\},m\in\{1,\ldots,k-1\},s\in[0,1/n])$ is
bounded in $\cg_{\val\alpha +4}(\R^d)$. Since $k\le \nt$ when
$(t^n_k,t)\in\ct_1$, this implies that $(\rho^{n,-}_{t^n_k},n\ge
2,k\in\{1,\ldots,n-1\})$ is $O(1/n)$ in $\cg_{\val\alpha+2}(\R^d)$.

Let us now prove the same for $\rho^{n,+}$. After $\val\alpha$
integrations by parts and after setting $z=y-\zeta\sqrt{t-s}$, 
we get that $((p_n-p)*_{g,\alpha}p)_+$ is a weighted sum of terms of the
form $\ci(\chi^{n,+}_{\alpha_1,\alpha_2})$ - see Lemma
\ref{lem:convol_esp} - with
\begin{multline*}
 \chi^{n,+}_{\alpha_1,\alpha_2}(s,t,x,y,\zeta) = \ind_{\ct_1^+}(s,t)(t-s)^{d/2}\partial^{\alpha_1}g(y-\zeta\sqrt{t-s})\\ \partial^{\alpha_2}_3(p_n-p)(s,x,y-\zeta\sqrt{t-s})p(t-s,y-\zeta\sqrt{t-s},y)
\end{multline*}
and $\val{\alpha_1}+\val{\alpha_2}=\val{\alpha}$. Now, from Corollary
\ref{cor:pn-p} in the appendix, $(p_n-p,n\ge 1)$ is $O(1/n)$ in $\cg_2(\R^d)$ so that, using the
same arguments as in Step 2 of the proof of Proposition \ref{prop:convol_esp}-(i), we
get that $((p_n-p)*_{g,\alpha}p)_+$ is $O(1/n)$ in
$\ch_{\val\alpha+2}(\R^d)$. Since $\rho^{n,+}_{t^n_k}(t,x,y) = \ind_{\ct_1^+}(t^n_k,t)
((p_n-p)*_{g,\alpha}p)_+(t^n_k,t,x,y)$, we conclude that $(\rho^{n,+}_{t^n_k},n\ge
2,k\in\{1,\ldots,n-1\})$ is $O(1/n)$ in $\cg_{\val\alpha+2}(\R^d)$.
\end{proof}

Lastly, starting from (\ref{eq:An3t}), Lemmas \ref{lem:densite_reste} and \ref{lem:densite_Psi} with $j=2$
imply

\begin{lemma}\label{lem:R^n_3,t}
Under (B) and (C), there exists a $O(1/n^2)$ sequence $(\pi^n_3,n\ge 1)$ in $\cg_4(\R^d)$ such that for all $t\in(0,1]$, $f\in\cpol$ and $x\in\R^d$,
\begin{equation*}\label{eq:densite_Rn2t}
A^n_{3,t}f(x)=\int_{\R^d}f(y)\pi^n_3(t,x,y)\: dy.
\end{equation*}
\end{lemma}

Statement (iii) of Theorem \ref{thm:principal} is now proved: it follows
from (\ref{eq:delta_nt_dev_2}), (\ref{eq:Rnt}) and Lemmas \ref{lem:C_t},
\ref{lem:R^n_1,t}, \ref{lem:R^n_2,t} and \ref{lem:R^n_3,t}.

We now also have all the tools to prove Theorem
\ref{thm:principal_km}. Indeed, note that (\ref{eq:delta_nt_dev}) combined with Lemmas \ref{lem:densite_reste} and \ref{lem:densite_Psi} imply that we have an expansion of arbitrary order $j$ for $p_n-p$:
\begin{equation*}
(p_n-p)(t,\cdot,\cdot) = \sum_{i=2}^{j}\inv{i!n^i}\sum_{k=0}^{\nt-1}
\psi^{n,i}_{t^{n}_{k}}(t,\cdot,\cdot)+ r^{n,j}(t,\cdot,\cdot)+\sum_{i=2}^{j}\frac{\left(t-\ntn\right)^i}{i!}\psi^{n,i}_{\ntn}(t,\cdot,\cdot).
\end{equation*}
Since $(r^{n,j},n\ge 1)$ is $O(1/n^j)$ in $\cg_{2j}(\R^d)$ and
$(\psi^{n,j}_{t^n_k},n\ge 1,k\in\{0,\ldots,n\})$ is bounded in
$\cg_{2j}(\R^d)$, this gives (\ref{eq:dev_densite_tout_ordre}) with $(\pi_{n,i},n\ge
1)$ bounded in $\cg_{2i-2}(\R^d)$ and $(\pi'_{n,i},n\ge 1)$ and
$(\pi''_{n,i},n\ge 1)$ bounded in $\cg_{2i}(\R^d)$.

\section{Appendix}\label{sec:appendix}

\subsection{Kernel of  $R^{n,j}$}

Here we make explicit the kernel of the remainder
$R^{n,j}_t$, recall (\ref{eq:R_njt}):






\begin{lemma}\label{lem:densite_reste}
Under (B) and (C), for each $j\in\N^*$, there exists a $O(1/n^j)$ sequence
$(r^{n,j},n\ge 1)$ in $\cg_{2j}(\R^d)$
such that for all $t\in(0,1]$, $n\ge 1$, $f\in\cpol$ and $x\in\R^d$,
\begin{equation*}\label{eq:densite_Phi}
R^{n,j}_{t}f(x)=\int_{\R^d}f(y)r^{n,j}(t,x,y)\: dy.
\end{equation*}
\end{lemma}

\begin{proof}
From (\ref{eq:R_njt}) and (\ref{eq:R_nkjt}), $R^{n,j}_{t} =
R^{n,j}_{1,t} + R^{n,j}_{2,t}$ where
\begin{eqnarray*}
R^{n,j}_{1,t} &=&
\sum_{k=0}^{\nt-1}\int_{0}^{1/n}\int_{0}^{s_{1}}\cdots\int_{0}^{s_{j}}P^n_{t^n_k}\Phi^{n,j+1}_{s_{j+1},1/n}P_{t-t^n_{k+1}}\:
ds_{j+1}\cdots ds_{2} ds_{1},\\
R^{n,j}_{2,t} &=&
\int_{0}^{t-\ntn}\int_{0}^{s_{1}}\cdots\int_{0}^{s_{j}}P^n_{\ntn}\Phi^{n,j+1}_{s_{j+1},t-\ntn}\:
ds_{j+1}\cdots ds_{2} ds_{1}.
\end{eqnarray*}
Let us first deal with $R^{n,j}_{1,t}$. Using the fact that $k\ge 1$ for
the first equality, (\ref{eq:def_Phi_prime}) for the second one, the fact that $P_{1/n-s}P_{t-t^n_{k+1}}=P_{t-t^n_{k}-s}$ for
the third one, and (\ref{eq:deriv_Ptf}) and Fubini's theorem for the
last one, we have for all $f\in\cpol$, $x\in\R^d$, $t\in[0,1]$, $n\ge
1$, $k\in\{1,\ldots,\nt-1\}$ and $s\in(0,1/n)$,
\begin{eqnarray*}
&&P^n_{t^n_k}\Phi^{n,j}_{s,1/n}P_{t-t^n_{k+1}}f(x) \\ 
&&=\int_{\R^d}p_n(t^n_k,x,z_1)\Phi^{n,j}_{s,1/n}P_{t-t^n_{k+1}}f(z_1)\: dz_1 \\
&&= \sum_{1\le |\alpha|\le 2j}
\sum_{l=1}^{m_{j,\alpha}}\int_{\R^d}p_n(t^n_k,x,z_1)
g_{j,\alpha,l}(z_1)P^n_{s}(h_{j,\alpha,l}\partial^{\alpha}P_{1/n-s})P_{t-t^n_{k+1}}f(z_1)\:
dz_1 \\
&&= \sum_{1\le |\alpha|\le 2j}
\sum_{l=1}^{m_{j,\alpha}}\int_{\R^d}\int_{\R^d}p_n(t^n_k,x,z_1)
g_{j,\alpha,l}(z_1)p_n(s,z_1,z_2)h_{j,\alpha,l}(z_2)\partial^{\alpha}P_{t-t^n_{k}-s}f(z_2)\:
dz_2 dz_1 \\
&&= \int_{\R^d}f(y)\vphi^{n,j}_{t^n_k}(s,t,x,y)\: dy
\end{eqnarray*}
where $\vphi^{n,j}_{t^n_k}= \sum_{1\le |\alpha|\le 2j}
\sum_{l=1}^{m_{j,\alpha}} \vphi^{n,j}_{t^n_k,\alpha,l}$ with
\begin{multline*}
\vphi^{n,j}_{t^n_k,\alpha,l}(s,t,x,y) = \ind_{(0,\inv
  n)}(s)\ind_{[t^n_{k+1},1]}(t) \\
\int_{\R^d}(p_n*_{g_{j,\alpha,l},0}p_n)(t^n_k,t^n_k+s,x,z_2)h_{j,\alpha,l}(z_2)\partial^{\alpha}_2p(t-t^n_k-s,z_2,y)\: dz_2.
\end{multline*}
Now, setting $q^{n,j}_{t^n_k,\alpha,l}(u,x,z)=\ind_{(t^n_{k},1]}(u)(p_n*_{g_{j,\alpha,l},0}p_n)(t^n_k,u,x,z)$, it follows from
Proposition \ref{prop:convol_esp}-(ii) that $(q^{n,j}_{t^n_k,\alpha,l},n\ge
1,k\in\{1,\ldots,n\})$ is a bounded sequence in $\cg(\R^d)$. Since
$\vphi^{n,j}_{t^n_k,\alpha,l}(s,t,x,y)= \ind_{(0,\inv n)}(s)\ind_{[t^n_{k+1},1]}(t)(q^{n,j}_{t^n_k,\alpha,l}*_{h_{j,\alpha,l},\alpha}p)(t^n_k+s,t,x,y)$,
Proposition \ref{prop:convol_esp}-(i) shows that $(\vphi^{n,j}_{t^n_k,\alpha,l},n\ge
1,k\in\{1,\ldots,n\})$ is bounded in $\ch_{\val\alpha}(\R^d)$, so that $(\vphi^{n,j}_{t^n_k},n\ge
1,k\in\{1,\ldots,n\})$ is bounded in $\ch_{2j}(\R^d)$.

When $k=0$, we have in the same way for all $f\in\cpol$
\begin{equation*}
\Phi^{n,j}_{s,1/n}P_{t-1/n}f(x) =  \int_{\R^d}f(y)\vphi^{n,j}_{0}(s,t,x,y)\: dy
\end{equation*}
where $\vphi^{n,j}_{0}= \sum_{1\le |\alpha|\le 2j}
\sum_{l=1}^{m_{j,\alpha}} \vphi^{n,j}_{0,\alpha,l}$ with
\begin{equation*}
\vphi^{n,j}_{0,\alpha,l}(s,t,x,y) = \ind_{(0,\inv n)}(s)\ind_{[\inv n,1]}(t)g_{j,\alpha,l}(x)(p_n*_{h_{j,\alpha,l},\alpha}p)(s,t,x,y).
\end{equation*}
Again Proposition \ref{prop:convol_esp}-(i) imply that $(\vphi^{n,j}_{0},n\ge
1)$ is bounded in $\ch_{2j}(\R^d)$.

Eventually, for all $f\in\cpol$, we have
$R^{n,j}_{1,t}f(x)=\int_{\R^d}f(y)r^{n,j}_1(t,x,y)\: dy$ with
\begin{equation*}
r^{n,j}_{1}(t,x,y) = \sum_{k=0}^{\nt-1}\int_{0}^{1/n}\int_{0}^{s_{1}}\cdots\int_{0}^{s_{j}}\vphi^{n,j+1}_{t^n_k}(s_{j+1},t,x,y)\:
ds_{j+1}\cdots ds_{2} ds_{1},
\end{equation*}
and since the family $(\vphi^{n,j+1}_{t^n_k},n\ge
1,k\in\{0,\ldots,n\})$ is bounded in $\ch_{2j+2}(\R^d)$, the sequence $(r^{n,j}_{1},n\ge 1)$ is $O(1/n^j)$ in $\cg_{2j}(\R^d)$.

As for $R^{n,j}_{2,t}$, similar arguments lead to
\begin{equation*}
P^n_{\ntn}\Phi^{n,j}_{s,t-\ntn}f(x) =  \int_{\R^d}f(y)\phi^{n,j}(s,t,x,y)\: dy
\end{equation*}
where $\phi^{n,j}= \sum_{1\le |\alpha|\le 2j}
\sum_{l=1}^{m_{j,\alpha}} \phi^{n,j}_{\alpha,l}$ with
\begin{multline*}
\phi^{n,j}_{\alpha,l}(s,t,x,y) = \ind_{[\inv
  n,1]}(t)\ind_{\left(0,t-\frac{\nt}{n}\right)}(s)\int_{\R^d}(p_n*_{g_{j,\alpha,l},0}p_n)\left(\frac{\nt}{n},\frac{\nt}{n}+s,x,z_2\right) \\ h_{j,\alpha,l}(z_2)\partial^{\alpha}_2p\left(t-\frac{\nt}{n}-s,z_2,y\right)\: dz_2 + \ind_{\{0<s<t<\inv n\}}g_{j,\alpha,l}(x)
\left(p_n*_{h_{j,\alpha,l},\alpha}p \right)(s,t,x,y).
\end{multline*}
We can treat $\phi^{n,j}_{\alpha,l}$ exactly as we have treated
$\vphi^{n,j}_{t^n_k,\alpha,l}$, and get that $(\phi^{n,j},n\ge 1)$ is
bounded in $\ch_{2j}(\R^d)$, so that $R^{n,j}_{2,t}$ has a kernel
$(r^{n,j}_{2},n\ge 1)$ defined by
\begin{equation*}
r^{n,j}_{2}(t,x,y) = \int_{0}^{t-\ntn}\int_{0}^{s_{1}}\cdots\int_{0}^{s_{j}}\phi^{n,j+1}(s_{j+1},t,x,y)\:
ds_{j+1}\cdots ds_{2} ds_{1}
\end{equation*}
which is $O(1/n^{j})$ in $\cg_{2j}(\R^d)$.

Eventually, putting $r^{n,j}=r^{n,j}_1+r^{n,j}_2$ completes the proof.
\end{proof}

In particular we have

\begin{corollary}\label{cor:pn-p}
Under (B) and (C), $(p_n-p,n\ge 1)$ is $O(1/n)$ in $\cg_2(\R^d)$.
\end{corollary}

\begin{proof}
From (\ref{eq:delta_nt_dev}) applied with $j=1$ and Lemma \ref{lem:densite_reste}, we have for all
$f\in\cpol$
\begin{equation*}
\int_{\R^d}f(y)(p_n-p)(t,x,y)\: dy=\Delta^n_tf(x)=R^{n,1}_tf(x)=\int_{\R^d}f(y)r^{n,1}(t,x,y)\: dy
\end{equation*}
so that $p_n-p=r^{n,1}$, and Lemma \ref{lem:densite_reste} gives the result.
\end{proof}

Eventually, we have kernels for the operators
$P^n_{t^{n}_{k}}L^*_jP_{t-t^{n}_{k}}$:

\begin{lemma}\label{lem:densite_Psi}
Under (B) and (C), for each $j\in\N^*$, there exists a bounded sequence
$(\psi^{n,j}_{t^n_k},n\ge 1,k\in\{0,\ldots,n\})$ in $\cg_{2j}(\R^d)$
such that for all $t\in(0,1]$, $n\ge 1$, $k\in\{0,\ldots,\nt\}$, $f\in\cpol$ and $x\in\R^d$,
\begin{equation*}\label{eq:densite_Psi}
P^n_{t^{n}_{k}}L^*_jP_{t-t^{n}_{k}}f(x)=\int_{\R^d}f(y)\psi^{n,j}_{t^n_k}(t,x,y)\: dy.
\end{equation*}
\end{lemma}

The proof is omitted since it copies the arguments of the proof of Lemma
\ref{lem:densite_reste} - it is even a bit simpler.

\subsection{Operators on $\cg_l(\R^d)$}\label{subsec:op_cgl}

When $\pi\in\cg_l(\R^d)$, $\pi(t,\cdot,y)\in L^{\infty}(\R^d)$ so that for $s\in[0,1]$ and $n\ge 1$ we can define two functions $P_s\pi$ and $P^n_s\pi$ on $(0,1]\times\R^d\times\R^d$ by $P_s\pi(t,\cdot,y)=\ind_{\{s\le t\}}P_s(\pi(t,\cdot,y))$ and $P^n_s\pi(t,\cdot,y)=\ind_{\{s\le t\}}P^n_s(\pi(t,\cdot,y))$, i.e.
\begin{equation}\label{eq:def_P_cgl}
P_s\pi(t,x,y)=\ind_{\{s\le t\}}\E\left[\pi\left(t,X^x_s,y\right) \right] \quad \mbox{and} \quad P^n_s\pi(t,x,y)=\ind_{\{s\le t\}}\E\left[\pi\left(t,X^{n,x}_s,y\right) \right].
\end{equation}
We also write $\Delta^n_s\pi=P^n_s\pi - P_s\pi$. For $j\in\N^*$ we denote by $\Phi^j$ the family $(\Phi^{n,j}_{s,1/n},n\ge 1,s\in[0,1/n])$ of operators on $\cg_l(\R^d)$ defined as in (\ref{eq:def_Phi}) by
\begin{equation*}
\Phi^{n,j}_{s,1/n}\pi(t,x,y)=\E\left[L^{x}_jP_{1/n-s}\pi\left(t,X^{n,x}_{s},y\right)\right],
\end{equation*}
i.e., using (\ref{eq:lxj}),
\begin{equation}\label{eq:def_Phi_2}
\Phi^{n,j}_{s,1/n}=\sum_{1\le |\alpha|\le 2j} \sum_{l=1}^{m_{j,\alpha}} g_{j,\alpha,l}P^n_{s}\left(h_{j,\alpha,l}\partial^{\alpha}P_{1/n-s}\right).
\end{equation}
Denoting by $\cl_b(\cg_l(\R^d),\cg_{l'}(\R^d))$ the space of all morphisms mapping any bounded subset of $\cg_{l}(\R^d)$ into a bounded subset of $\cg_{l'}(\R^d)$, we then have

\begin{proposition}\label{prop:operators_on_cgl}
Under (B) and (C), $(P_s,s\in[0,1])$ and $(P^n_s,s\in[0,1],n\ge1)$ are bounded families in $\cl_b(\cg_l(\R^d))$, and $\Phi^j$ is a bounded family in $\cl_b(\cg_l(\R^d),\cg_{l+2j}(\R^d))$.
\end{proposition}

\begin{proof}
Let us first deal with $(P_s)$. Let $\pi\in\cg_l(\R^d)$. $P_s$ is measurable. Moreover, Lebesgue's dominated convergence theorem shows that $P_s\pi(t,x,\cdot)$ is infinitely differentiable and that for all $\beta\in\N^d$
\begin{equation*}
\partial^{\beta}_yP_s\pi(t,x,y)=\ind_{\{s\le t\}}\E\left[\partial^{\beta}_3\pi\left(t,X^x_s,y\right) \right].
\end{equation*}
Hypothesis (A) ensures that a version of $X^x$ can be chosen such that for
each $t\ge 0$, the map $x\mapsto X^x_t$ is infinitely differentiable
(see, for example, \cite{kunita}). Since $\partial^{\beta}_3\pi(t,\cdot,y)\in\cpol$, it follows from Theorem 3.14 page 16 in \cite{kusuoka-stroock} that $\partial^{\beta}_yP_s\pi(t,\cdot,y)$ is infinitely differentiable and that for all $\alpha\in\N^d$ there exists universal polynomials $(\Pi_{\alpha,\mu},\val\mu\le\val\alpha)$ such that 
\begin{equation}\label{eq:demo_-1_operators_on_cgl}
\partial^{\alpha}_x\partial^{\beta}_yP_s\pi(t,x,y)=\ind_{\{s\le t\}}\sum_{\val\mu\le\val\alpha}\E\left[\partial^{\mu}_2\partial^{\beta}_3\pi\left(t,X^x_s,y\right)\Pi_{\alpha,\mu}\left(\partial^{\nu}_xX^x_s,\val\nu\le\val\alpha \right)\right]
\end{equation}
with
\begin{equation}\label{eq:demo_-0.5_operators_on_cgl}
\sup_{s\in[0,1],x\in\R^d}\E[\Pi_{\alpha,\mu}\left(\partial^{\nu}_xX^x_s,\val\nu\le\val\alpha \right)^2]<\infty
\end{equation}
 for all $\val\mu\le\val\alpha$. As a consequence, $P_s\pi(t,\cdot,\cdot)$ is infinitely differentiable and using Cauchy-Schwarz's inequality, (\ref{eq:def_cgl1}) and (\ref{eq:demo_-0.5_operators_on_cgl}), we see that for all bounded $\cb\subset\cg_l(R^d)$ and $\alpha,\beta\in\N^{d}$, there exists two constants $c_1\ge 0$ and $c_2>0$ such that for all $\pi\in\cb$, $s\in[0,1]$, $t\in(0,1]$ and $x,y\in\R^d$,
\begin{equation}\label{eq:demo_0_operators_on_cgl}
\val{\partial^{\alpha}_x\partial^{\beta}_yP_s\pi(t,x,y)} \le c_1 \ind_{\{s\le t\}}t^{-(\val\alpha+\val\beta+d+l)/2}\E\left[\exp\left(-c_2\norm{X^s_x-y}^2/t\right)\right]^{1/2}.
\end{equation}
Now, partitioning $\Omega$ into $\{\norm{X^s_x-y}\le\norm{x-y}/2\}$ and
$\{\norm{X^s_x-y}>\norm{x-y}/2\}$, we have
\begin{equation}\label{eq:demo_1_operators_on_cgl}
\E\left[\exp\left(-c_2\norm{X^s_x-y}^2/t\right)\right] \le \P\left(\norm{X^s_x-y}\le\norm{x-y}/2\right) + \exp\left(-c_2\norm{x-y}^2/4t\right).
\end{equation}
Using
(\ref{eq:def_cg2}) for $p\in\cg(\R^d)$ for
the fourth inequality, we can find $c_3,c_5\ge 0$ and $c_4,c_6>0$ such that for all
$s\in(0,1]$ and $x,y\in\R^d$,
\begin{eqnarray}
\P\left(\norm{X^x_s-y}\le\norm{x-y}/2\right) &\le&
\P\left(\norm{X^x_s-x}\ge\norm{x-y}/2\right) \nonumber\\
&=& \int_{\R^d}\ind_{\{\norm{z-x}\ge\norm{x-y}/2\}}p(s,x,z)\: dz\nonumber\\
&=& \int_{\R^d}\ind_{\{\norm
  \xi\ge\norm{x-y}/2\sqrt{s}\}}p(s,x,x+\xi\sqrt{s})s^{d/2}\: d\xi\nonumber\\
&\le& c_3\int_{\R^d}\ind_{\{\norm
  \xi\ge\norm{x-y}/2\sqrt{s}\}}\exp(-c_4\norm{\xi}^2)\: d\xi\nonumber\\
&\le& c_5\exp\left(-c_6\norm{x-y}^2/s\right). \label{eq:demo_2_operators_on_cgl}
\end{eqnarray}
Eventually, from (\ref{eq:demo_1_operators_on_cgl}) and (\ref{eq:demo_2_operators_on_cgl}), we can find $c_7\ge 0$ and $c_8>0$ such that for all $s\in[0,1]$, $t\in(0,1]$ and $x,y\in\R^d$,
\begin{eqnarray}
\ind_{\{s\le t\}}\E\left[\exp\left(-c_2\norm{X^s_x-y}^2/t\right)\right] &\le&
c_5\exp\left(-c_6\norm{x-y}^2/t\right)+
\exp\left(-c_2\norm{x-y}^2/4t\right) \nonumber\\
&\le& c_7\exp\left(-c_8\norm{x-y}^2/t\right). \label{eq:demo_3_operators_on_cgl}
\end{eqnarray}
It is enough to inject (\ref{eq:demo_3_operators_on_cgl}) into
(\ref{eq:demo_0_operators_on_cgl}) to complete the proof for $(P_s)$.

This proof naturally extends to the case of $(P^n_s)$. Indeed,
(\ref{eq:demo_-1_operators_on_cgl}) holds with $(X^n,P^n)$ instead
of $(X,P)$. Moreover, from Lemma \ref{lem:moments_euler_B}, (\ref{eq:demo_-0.5_operators_on_cgl}) holds uniformly in $n$ with $X^n$ instead of $X$. Eventually, (\ref{eq:demo_2_operators_on_cgl}) holds with $X^n$ instead of $X$, uniformly in $n$ because $(p_n,n\ge1)$ is bounded in $\cg(\R^d)$.

As for $\Phi^j$, it is enough to use (\ref{eq:def_Phi_2}), the
boundedness of $(P_s)$ and $(P^n_s)$, Remark \ref{rem:cbo} and the facts that multiplication by a function in $\cb$ belongs to $\cl_b(\cg_l(\R^d),\cg_{l}(\R^d))$ and that $\partial^{\alpha}_2\in\cl_b(\cg_l(\R^d),\cg_{l+\val\alpha}(\R^d))$.
\end{proof}

\subsection{Moments for the Euler scheme and its derivatives}\label{subsec:moments_euler}

Let us assume (A). Then it is known that $X^{n,x}_t$ has bounded moments of any order and that for all $q\in\N$, one can find $c\ge 0$ such that for all $x\in\R^d$,
\begin{equation}\label{eq:moments_euler}
\sup_{t\in[0,1],n\ge 1}\E\left[\norm{X^{n,x}_t}^q \right]\le c\left(1+\norm x ^q\right)
\end{equation}
(see \cite{talay2}). From (\ref{eq:def_euler}), $x\mapsto X^{n,x}_t$ is infinitely differentiable and we shall see that analogous upper bounds hold for its derivatives. Following \cite{kusuoka-stroock}, for $m\ge 1$, we denote by $X^{(m),n,x}_t$ the $m$-th derivative of $x\mapsto X^{n,x}_t$ at point $x$. It should be thought of as a $d\times d^m$ matrix. For instance, $X^{(1),n,x}_t$ is the jacobian matrix of $x\mapsto X^{n,x}_t$. Differentiating (\ref{eq:def_euler}), we have
\begin{equation}\label{eq:EDSderiv_1}
X^{(1),n,x}_t=I+\int_{0}^{t}b^{(1)}(X^{n,x}_{\nsn})X^{(1),n,x}_{\nsn}\:ds+\sum_{j=1}^r\int_{0}^{t}\sigma_j^{(1)}(X^{n,x}_{\nsn})X^{(1),n,x}_{\nsn}\:dB^j_{s},
\end{equation}
where $I$ stands for the identity matrix and $\sigma_j$ is the $j$-th column of $\sigma$. Besides,
by induction, there are for each $m\ge 2$ universal polynomials $P_{m,j}$, $j\in\{0,\ldots,r\}$, such that
 \begin{multline}\label{eq:EDSderiv_m}
X^{(m),n,x}_t=\int_{0}^{t}b^{(1)}(X^{n,x}_{\nsn})X^{(m),n,x}_{\nsn}\:ds+\sum_{j=1}^r\int_{0}^{t}\sigma_j^{(1)}(X^{n,x}_{\nsn})X^{(m),n,x}_{\nsn}\:dB^j_{s} \\ + \int_{0}^{t}Q^{n,x}_{m,0,\nsn}\: ds + \sum_{j=1}^r\int_{0}^{t}Q^{n,x}_{m,j,\nsn}\: dB^j_s,
\end{multline}
where
\begin{eqnarray}\label{eq:polQ}
\left\{
          \begin{array}{rll}
Q^{n,x}_{m,0,t}&=&P_{m,0}(b^{(2)}(X^{n,x}_t),\ldots,b^{(m)}(X^{n,x}_t),X^{(1),n,x}_t,\ldots,X^{(m-1),n,x}_t),\\
Q^{n,x}_{m,j,t}&=&P_{m,j}(\sigma_j^{(2)}(X^{n,x}_t),\ldots,\sigma_j^{(m)}(X^{n,x}_t),X^{(1),n,x}_t,\ldots,X^{(m-1),n,x}_t).
          \end{array}
        \right.
\end{eqnarray}
This is analogous to (1.8) page 4 in \cite{kusuoka-stroock}. Then we have

\begin{lemma}\label{lem:moments_euler_A}
Under (A), for all $m\ge 1$ and $q\in\N$, there exists $c\ge 0$ and $q'\in\N$ such that for all $x\in\R^d$,
\begin{equation}\label{eq:mom_euler_A}
\sup_{t\in[0,1],n\ge 1}\E\left[\norm{X^{(m),n,x}_t}^q \right]\le c\left(1+\norm x ^{q'}\right).
\end{equation}
\end{lemma}

\begin{proof}
We give a proof by induction on $m$. Let us first assume that $m=1$. Let $q\in\N$. From (\ref{eq:EDSderiv_1}), and observing that (A) states that $b^{(1)}$ and all the $\sigma_j^{(1)}$ are bounded, Jensen's and Burkholder-Davis-Gundy's inequalities lead to the existence of $c\ge 0$ such that for all $t\in[0,1]$, $n\ge 1$ and $x\in\R^d$,
\begin{equation*}
\E\left[\norm{X^{(1),n,x}_t}^q\right]\le c \left(1+\int_{0}^{t}\E\left[\norm{X^{(1),n,x}_{\nsn}}^q\right]\: ds\right).
\end{equation*}
Taking this inequality at time $\ntn$ and applying Gronwall's lemma, we get that
\begin{equation*}
\sup_{t\in[0,1],n\ge 1,x\in\R^d}\E\left[\norm{X^{(1),n,x}_{\ntn}}^q \right]<\infty.
\end{equation*}
From (\ref{eq:def_euler}), one easily checks that the same holds at time $t$ instead of $\ntn$, so that (\ref{eq:mom_euler_A}) holds for $m=1$ with $q'=0$.

Let us now assume that (\ref{eq:mom_euler_A}) holds for the $m-1$ first derivatives. Let $q\in\N$. From (\ref{eq:EDSderiv_m}), and observing again that (A) states that $b^{(1)}$ and all the $\sigma_j^{(1)}$ are bounded, Jensen's and Burkholder-Davis-Gundy's inequalities lead to the existence of $c_1\ge 0$ such that for all $t\in[0,1]$, $n\ge 1$ and $x\in\R^d$,
\begin{equation}\label{eq:moments_euler_A_proof}
\E\left[\norm{X^{(m),n,x}_t}^q\right]\le c_1 \left(\int_{0}^{t}\E\left[\norm{X^{(m),n,x}_{\nsn}}^q\right]\: ds+\int_{0}^{t}\sum_{j=0}^r\E\left[\norm{Q^{n,x}_{m,j,\nsn}}^q\right]\: ds\right).
\end{equation}
Using (\ref{eq:polQ}), the induction hypothesis, (A) and (\ref{eq:moments_euler}), we find $c_2\ge 0$ and $q'\in\N$ such that for all $s\in[0,1]$, $n\ge 1$ and $x\in\R^d$,
\begin{equation*}
\sum_{j=0}^r\E\left[\norm{Q^{n,x}_{m,j,\nsn}}^q\right] \le c_2\left(1+\norm x ^{q'} \right).
\end{equation*}
Thus, taking (\ref{eq:moments_euler_A_proof}) at time $\ntn$ and applying Gronwall's lemma, we find $c\ge 0$ such that for all $x\in\R^d$,
\begin{equation*}
\sup_{t\in[0,1],n\ge 1}\E\left[\norm{X^{(m),n,x}_{\ntn}}^q \right]\le c\left(1+\norm x ^{q'}\right).
\end{equation*}
From (\ref{eq:def_euler}), one easily checks that the same holds at time $t$ instead of $\ntn$, which completes the proof.
\end{proof}

Observe that, under (B), the above proof holds with $q'=0$ so that we have

\begin{lemma}\label{lem:moments_euler_B}
Under (B), for all $m\ge 1$ and $q\in\N$,
\begin{equation*}\label{eq:moments_euler_A}
\sup_{t\in[0,1],n\ge 1,x\in\R^d}\E\left[\norm{X^{(m),n,x}_t}^q \right]<\infty.
\end{equation*}
\end{lemma}

{\bf Acknowledgement.} The author would like to thank the referee and
J.-F. \textsc{Delmas} for
thorough readings of the paper and for many useful suggestions for
improvement. He also wishes to thank V. \textsc{Bally} for his encouragement.




\newcommand{\sortnoop}[1]{}

\end{document}